\documentclass[12pt]{article}
\usepackage{amsmath}
\usepackage{amsfonts}
\usepackage{graphicx,psfrag,epsf}
\usepackage{enumerate}
\usepackage[round]{natbib}
\usepackage{url} 
\usepackage{amssymb}
\usepackage{bbm}
\usepackage{caption}
\usepackage{subcaption}
\usepackage{siunitx}
\usepackage{float}
\usepackage{ulem} 
\usepackage{comment} 
\usepackage{adjustbox} 

\newcommand{\blind}{0}
\usepackage[usenames, dvipsnames]{xcolor}
\newcommand{\given}{\,|\,}
\NewDocumentCommand{\evalat}{sO{\big}mm}{%
  \IfBooleanTF{#1}
   {\mleft. #3 \mright|_{#4}}
   {#3#2|_{#4}}%
}

\usepackage{amsthm}
\newtheorem{theorem}{Theorem}[section]
\newtheorem{corollary}{Corollary}[theorem]
\newtheorem{lemma}[theorem]{Lemma}
\newtheorem{proposition}[theorem]{Proposition}

\newcommand{\LZ}[1]{\textcolor{RedViolet}{[LZ: #1]}}
\newcommand{\edit}[1]{\color{black}{#1}}
\newcommand{\editLZ}[1]{\color{Cerulean}{#1}\color{black}} 
\def\BASadd#1{{\textcolor{blue}{(BAS: #1)}}}

\def\replace#1{{\textcolor{black}{#1}}}

\renewcommand{\P}{\mbox{P}}
\renewcommand{\Pr}{\mbox{P}}

\newcommand{\RV}{\mbox{RV}}

\DeclareMathOperator{\E}{E}
\newcommand{\btheta}{ {\boldsymbol \theta} }

\usepackage{hyperref}
\hypersetup{
    colorlinks=true,
    linkcolor=blue,
    filecolor=magenta,      
    urlcolor=blue,
}

\begin{document}
\def\spacingset#1{\renewcommand{\baselinestretch}%
{#1}\small\normalsize} \spacingset{1}

\hypersetup{citecolor=blue}
\if0\blind
{
  \title{\bf Hierarchical Transformed Scale Mixtures for Flexible Modeling of Spatial Extremes on Datasets with Many Locations}
  \author{Likun Zhang\\
    Department of Statistics,  Pennsylvania State University\\
    Benjamin A. Shaby\\
    Department of Statistics,  Colorado State University\\
    Jennifer L. Wadsworth \\
    Department of Mathematics and Statistics, Fylde College, Lancaster University
    }
  \maketitle
} \fi

\if1\blind
{
  \bigskip
  \bigskip
  \bigskip
  
    \title{\bf Hierarchical Transformed Scale Mixtures for Flexible Modeling of Spatial Extremes on Datasets with Many Locations}
    \maketitle

  \medskip
} \fi

\bigskip
\begin{abstract}
Flexible spatial models that allow transitions between tail dependence classes have recently appeared in the literature. However, inference for these models is computationally prohibitive, even in moderate dimensions, due to the necessity of repeatedly evaluating the multivariate Gaussian distribution function. In this work, we attempt to achieve truly high-dimensional inference for extremes of spatial processes, while retaining the desirable flexibility in the tail dependence structure, by modifying an established class of models based on scale mixtures Gaussian processes. We show that the desired extremal dependence properties from the original models are preserved under the modification, and demonstrate that the corresponding Bayesian hierarchical model does not involve the expensive computation of the multivariate Gaussian distribution function. We fit our model to exceedances of a high threshold, and perform coverage analyses and cross-model checks to validate its ability to capture different types of tail characteristics. We use a standard adaptive Metropolis algorithm for model fitting, and further accelerate the computation via parallelization and \texttt{Rcpp}. Lastly, we apply the model to a dataset of a fire threat index on the Great Plains region of the US, which is vulnerable to massively destructive wildfires. We find that the joint tail of the fire threat index exhibits a decaying dependence structure that cannot be captured by limiting extreme value models.
\end{abstract}

\vfill

\newpage
\spacingset{1.45} 

\section{Introduction}\label{chap1:intro}
Modeling the dependence structure in the extremes of spatial processes is of great consequence for risk analysis of extreme events. In this paper, we make a slight alteration to the flexible class of randomly scaled transformed Gaussian process models to sidestep computational bottlenecks normally encountered in the likelihood. Our modification enables high-dimensional inference, while preserving submodels that transition smoothly between extremal dependence classes.

Generally, the probability that two spatially-indexed random variables exceed a high level simultaneously varies by their separation distance, and the particular way in which this joint probability decays must be well-represented in models if one hopes to accurately to assess risks posed by spatial extremal phenomena. Good estimation of how the dependence changes both with distance in space and as \replace{one moves farther into }the joint tail will enable us to accurately calculate exceedance probabilities of areal quantities, predict at un-observed locations, and, secondarily, get a more realistic picture of marginal quantities.

In classical spatial modeling, Gaussian processes have been widely used due to their mathematical simplicity and tractability for larger datasets. However, the Gaussian density function is very light-tailed, and thus \replace{has the potential }to underestimate probabilities associated with extreme events; furthermore, Gaussian models stipulate that the dependence among rare events at distinct locations will always diminish such that the probability of observing an extreme at one location, conditional upon an extreme at another location, is zero in the limit. This property is termed asymptotic independence, but models that only exhibit this phenomenon may be too inflexible for applications where the true tail dependence structure is uncertain. 

Max-stable processes form an important class of models that exhibit the alternative scenario of asymptotic dependence. They are the natural extension of classical univariate extreme value theory to infinite dimensional settings, and therefore can provide an asymptotically-justified modeling framework for datasets consisting of block-maxima. Counterparts of max-stable processes suitable for threshold exceedances are called generalized Pareto processes (\citealp{ferreira2014generalized}; \citealp{thibaud2015efficient}).  These processes are also asymptotically dependent, but possess the advantage that they bypass many of the computational difficulties of max-stable processes. 

Despite the theoretical appeal of limiting max-stable and generalized Pareto processes, there are two main drawbacks to these models: (i) the assumption of asymptotic dependence may be incorrect and (ii) even if the data are asymptotically dependent, they will often not appear to follow such limiting models at sub-asymptotic levels. Both max-stable and generalized Pareto process dependence structures exhibit stability properties, meaning that their dependence structures are invariant to the operations of taking maxima and conditioning upon exceedances of higher thresholds, respectively. If the true data generating process exhibits weakening dependence in the un-observed region of the tail, inference drawn under these models about the far joint tail will over-estimate risk, sometimes substantially.

On account of the limitations of limiting models, it is desirable to find a family of spatial models that can transition between asymptotic dependence and asymptotic independence. In particular, we will be examining a class of marginally transformed Gaussian scale mixture models, which includes those recently proposed by \citet{huser2017bridging} and \citet{huser2017modeling}. These models are of great interest due to their appealing theoretical properties and ability to flexibly capture both types of extremal dependence structure. Unfortunately, inference for these models is not feasible for \replace{large numbers of observed sites}, as calculation of the censored likelihood, the preferred method for fitting joint tail models, entails integration over high-dimensional multivariate Gaussian distribution functions. To increase scalability, we propose an adaptation of the model by adding an independent measurement error term to each component. By adding this nugget effect, the new model circumvents the lengthy computation of the multivariate normal distribution function. Also it can elegantly avoid the integral of the process below the censoring threshold by considering the uncensored process as latent and drawing from it using Gibbs sampling, allowing for truly high-dimensional inference. Furthermore, we show that the modified models retain all the significant asymptotic properties of the original smooth models despite the presence of the measurement errors, which lays a solid theoretical foundation for correctly capturing the sub-asymptotic dependence behavior.  

The article is organized as follows. Section \ref{chap1:lr} provides a brief literature review on the measures of extremal dependence and hybrid spatial extreme models, and further explains the intractability of the existing censored likelihood approaches for inference on these hybrid models. Section \ref{sec:model} describes our new model that alleviates the computational problems, and studies its extremal dependence properties. Section \ref{sec:bayes} includes a marginal transformation in the hierarchical model and details the inference using Gibbs sampling. Section \ref{sec:simu} presents a simulation study that validates the methodology. We apply our model to a dataset of the Fosberg Fire Weather Index (FFWI) on the Great Plains in Section \ref{data_analysis}. Section \ref{sec:discussion} concludes with some discussion. Appendix \ref{app:proofs} provides 
proofs of all the theoretical results. Appendix \ref{sec:diagnotics} includes supplementary diagnostics for the data application.

\section{Spatial Dependence for Extremes}\label{chap1:lr}
For a stochastic process $\{X(\boldsymbol{s}):\boldsymbol{s}\in \mathcal{S}\}$, we write $X_j=X(\boldsymbol{s}_j)$ and so forth for simplicity, where $\boldsymbol{s}_j$ denotes the $j$th spatial location. \replace{It is useful to summarize the extremal dependence implied by the observed process concisely.}

We restrict the scope to the bivariate case, focusing on stationary and isotropic random fields. One example of a bivariate dependence measure is the upper tail dependence coefficient:
\begin{equation}\label{chiu}
    \chi_u(h)=P(F_j(X_j)>u \given F_k(X_k)>u)
\end{equation}
where $X_j\sim F_j$, $X_k\sim F_k$, and $h=\|\boldsymbol{s}_j-\boldsymbol{s}_k\|$. \citet{joe1993parametric} defined the upper tail dependence parameter as the limit $\chi(h)=\lim_{u\rightarrow 1}\chi_u(h)$. Asymptotic dependence is attained if and only if $\chi(h)>0$, while $\chi(h)=0$ defines asymptotic independence.

For max-stable processes, \replace{$\chi_u(h)=2-V(1,1)+O(1-u), \;u \to 1$}, where $V(\cdot, \cdot) = \log F_{jk}(\cdot, \cdot)$. Max-stable distributions can be associated to a generalized Pareto counterpart, for which $\chi_u(h)\equiv \chi(h) = 2-V(1,1)$ for all $u$ above a certain level \citep{rootzen2018multivariate}. The fact that $\chi_u(h)$ does not depend on $u$ is a manifestation of the threshold-stability of generalized Pareto processes. In practice, empirical estimates of~\eqref{chiu} from data tend to show $\chi_u(h)$ decreasing with both $h$ and $u$, meaning that realistic models should also have this property.

When $\chi(h)=0$, i.e., the case of asymptotic independence, further detail about the behavior of $\chi_u(h)$ is obtained by exploiting the joint tail assumption of  \citet{ledford1996statistics}:
\begin{equation}\label{tailDependence}
    P(F_j(X_j)>u|F_k(X_k)>u)=\mathcal{L}(1-u)(1-u)^{1/\eta_X(h)-1}
\end{equation}
where $\mathcal{L}$ is slowly varying at zero, that is, $\lim_{t\rightarrow 0}\mathcal{L}(tx)/\mathcal{L}(t)=1$ for any $x>0$, and $\eta_X(h)\in (0,1]$ is the coefficient of tail dependence of the process $X$. The pair of variables $(X_j, X_k)$ are asymptotically dependent when $\eta_X(h)=1$ and $\mathcal{L}(\cdot)\nrightarrow 0$. The remaining cases are all asymptotically independent, and the value of $\eta_X(h)$ characterizes the strength of extremal dependence in the upper joint tail. In the case of a Gaussian process, $\eta_X(h)=\{1+\rho(h)\}/2$, where $\rho(h)$ is the correlation at lag $h$. The variables are called positively associated when $\eta_X(h)>1/2$ and negatively associated when $\eta_X(h)<1/2$. Near independence corresponds to $\eta_X(h)=1/2$.

Gaussian processes are asymptotically independent for all correlations $\rho(h)\neq 1$. They might be considered candidates for modeling the joint tail of asymptotically independent phenomena, but as there is no theory to specifically recommend Gaussian processes in this scenario, it is desirable to consider other models as well. As an alternative, \citet{opitz2016modeling} captures spatial dependence in asymptotically independent processes by construction of Laplace random fields, defined as mixtures of Gaussian processes with a random variance that is exponentially distributed. \citet{wadsworth2012dependence} proposed the class of inverted max-stable processes, for which the tail decay is specified fully by $\eta_X(h)$, although inference is computationally challenging.  

In real datasets, it is difficult to conclude definitively whether data exhibit asymptotic independence or asymptotic dependence, and incorrectly assuming an asymptotically independent model can lead to equally severe problems with bias as incorrectly assuming an asymptotically dependent model. Because of this, a recent focus in the literature has been on models that can encompass both scenarios.

\subsection{Traversing Asymptotic Independence and Dependence in Spatial Extremes}\label{hybrid}
\citet{wadsworth2012dependence} were the first to introduce hybrid models that combine max-stable and inverted max-stable processes so that asymptotic dependence prevails at short distances, and asymptotic independence at long distances. However, inference for this model is difficult because there are a fairly large number of parameters involved, and the transition between the dependence classes takes place at the boundary of the parameter space.

Recently, several Gaussian scale mixture models were proposed to allow more flexible transitions between dependence classes. Through multiplying an asymptotically independent Gaussian process by a random effect that governs the extremal dependence, these models can be described by a small number of parameters and have non-trivial asymptotically independent and asymptotically dependent submodels. More precisely, suppose $\{Z(\boldsymbol{s}),\boldsymbol{s}\in \mathcal{S}\}$ is a standard isotropic and stationary Gaussian process with covariance function 
$C_{\btheta_C}(h)$ indexed by a parameter vector $\btheta_C$, where $h$ is the length of the separation vector, so that $\Sigma_{\btheta_C}$ is the covariance matrix of associated finite-dimensional distributions. The class of Gaussian scale mixture models can be constructed as 
\begin{equation}\label{scalemixtureModel}
    X^*(\boldsymbol{s})=R\cdot g(Z(\boldsymbol{s})),\qquad R \given \btheta_R\sim F_R, 
\end{equation}
where $g(\cdot)$ is a link function, and $R>0$ is a random scaling factor, from distribution $F_R$ indexed by $\btheta_R$, that can be interpreted as a constant random process over spatial domain $\mathcal{S}$ with perfect dependence. Impacting simultaneously the whole domain $\mathcal{S}$, heavier tailed $R$ induces asymptotic dependence \edit{in $X^*$}, whereas lighter tailed $R$ induces asymptotic independence. \edit{\citet{engelke2018extremal} provide a fuller description of how extremal dependence of $X^*$ relates to the relative marginal tail heaviness of $R$ and $g(Z)$.}

\citet{morris2017space} uses a space-time model based on skew-$t$ process, where $g(\cdot)$ is a identity function, $R^2\sim \text{IG}(a/2,b/2)$ is an inverse gamma random variable, and $C_{\btheta_C}$ is a Mat\'ern covariance function. On top of the mixture, they added covariate effects and a skew term. Since the inverse gamma distribution is heavy tailed, the skew-$t$ process is asymptotically dependent for $a<\infty$. Asymptotic independence is achieved only when $a\rightarrow\infty$.

\citet{huser2017bridging} also used an identity link function, but placed few assumptions on the random scale, and provided more general results on the joint tail decay rates of the mixture processes. They showed that a wide class of Weibull-like tail decay in $R$ yields asymptotic independence, while a Pareto-like tail that is regularly varying at infinity gives asymptotic dependence. They also proposed a parametric model that bridges the two asymptotic regimes and provides a simple transition, in which $R$ is a two-parameter distribution
\begin{equation}
\label{hot}
    F_R(r)=\left\{
                \begin{array}{ll}
                  1-\exp\{-\gamma(r^{\beta}-1)/\beta\},&\beta>0,\\
                  1-r^{-\gamma},&\beta=0
                \end{array}
              \right.
\end{equation}
where $\gamma>0$, and the support is $[1,\infty)$. Since $(r^{\beta}-1)/\beta$ converges to $\log r$ as $\beta$ approaches 0, \eqref{hot} forms a continuous parametric family on $\beta$. When $\beta>0$, \eqref{hot} constitutes a class of Weibull-type distributions and thus assures asymptotic independence. When $\beta=0$, the variable $R$ is Pareto distributed and thus gives asymptotic dependence. This shows that the model provides greater flexibility and can transition from asymptotic dependence to independence via adjusting the value of $\beta$. 

However, the previous two Gaussian scale mixture models both make the transition between the dependence classes at the limit or the boundary of the parameter space. They are also inflexible in their representation of asymptotic dependence structures \replace{because there is dominating preference over one dependence class}. It may be more desirable to find a model for which the transition takes place in the interior of the parameter space so one we can quantify the uncertainty about the dependence class in a simpler manner. To overcome this, \citet{huser2017modeling} proposed a marginally transformed Gaussian scale mixture model, where $g(\cdot)$ transforms a standard Gaussian variable to standard Pareto, and $R$ itself is Pareto distributed:
\begin{equation}\label{huserWadsworth}
    g(z)=\frac{1}{1-\Phi(z)},\; R \given \delta \sim \text{Pareto}\left(\frac{1-\delta}{\delta}\right),\;\delta\in [0,1] .
\end{equation}
Here the type of asymptotic dependence is determined by the value of $\delta$. When $\delta\leq 1/2$, $R$ is light\replace{er tailed or equivalent to standard Pareto}, which induces asymptotic independence; when $\delta>1/2$, the converse is true, which induces asymptotic dependence. Specifically, the upper tail dependence parameter $\chi_{X^*}=\frac{2\delta-1}{\delta}E\left[\min\{g(Z_i),g(Z_k)\}^{(1-\delta)/\delta}\right]$ when $\delta>1/2$ and 0 otherwise, while the coefficient of tail dependence is
\[
   \eta_{X^*}=\left\{
                \begin{array}{ll}
                 1,&\delta>\frac{1}{2},\\
                 \frac{\delta}{1-\delta},&\frac{\eta_{Z}}{\eta_{Z}+1}<\delta\leq\frac{1}{2},\\
                 \eta_{Z},&\delta\leq\frac{\eta_{Z}}{\eta_{Z}+1},
                \end{array}
              \right.
  \]
where $\eta_{Z}$ is the coefficient of tail dependence for $(Z_i,Z_k)$ \citep{huser2017modeling}.

The model in \eqref{huserWadsworth} provides a smooth transition through asymptotically independent and asymptotically dependent submodels. It has many appealing asymptotic properties. However, inference for models of the form \eqref{scalemixtureModel} is typically made via censored likelihood. This requires computing an integral where the integrand contains the Gaussian distribution function in $|\mathcal{C}|$ dimensions, where $|\mathcal{C}|$ is the number of components below a designated high threshold. Such integrals are computationally prohibitive for even moderately-sized datasets. In Section~\ref{sec:model} we introduce a slight alteration to this model to make it tractable while preserving all the desired asymptotic results.

\subsection{The Censored Likelihood}\label{censored}
In multivariate and spatial extremes, the preferred approach to fitting the dependence structure is using a censored likelihood, which prevents observations from the bulk of the distribution from affecting the estimation of the extremal dependence structure. It provides a reasonable compromise between bias and variance compared to alternative approaches, although different censoring schemes have been adopted \citep{thibaud2015efficient,huser2016likelihood}.

For a process of the form \eqref{scalemixtureModel} observed at $D$ spatial locations $\boldsymbol{s}_1,\cdots, \boldsymbol{s}_D\in \mathcal{S}$, we obtain the distribution function by conditioning on $R$ as
\begin{equation}\label{cdf}
    G(\boldsymbol{x}^*)=\int_1^{r^*}\Phi_D\left(g^{-1}\left(\frac{\boldsymbol{x}^*}{r}\right);\mathbf{\Sigma}_{\theta_C}\right)f_R(r)dr,
\end{equation}
where $r^*=\min(x^*_1,\cdots, x^*_D)$, and $\Phi_D$ denotes the $D$-variate Gaussian distribution with zero mean and covariance matrix $\mathbf{\Sigma}_{\btheta_C}$.

Let $\mathcal{C}\subseteq \{1,\ldots, D\}$ be the set of locations with censored observations---that is, the set of locations where the components are below a high threshold; let $\mathcal{U}$ be the set of locations with uncensored observations. For any index set $A,B\subset \{1,\ldots, D\}$, denote $\boldsymbol{x}_A=\{\boldsymbol{x}_i:i\in A\}$, $\mathbf{\Sigma}_{A;B}$ as the matrix $\mathbf{\Sigma}$ restricted to the rows in $A$ and the columns in $B$, and let $\mathbf{\Sigma}_{A|B}$ be the Schur complement of $B$ in $\mathbf{\Sigma}_{A;B}$. The likelihood is obtained via taking partial derivatives of \eqref{cdf} with respect to $\mathcal{U}$:
\begin{equation}\label{censorLik}
    \begin{split}
       \frac{\partial ^{|\mathcal{U}|}}{\partial \boldsymbol{x}^*_{\mathcal{U}}}G(\boldsymbol{x}^*)
        &=\int_1^{r^*}\Phi_{|\mathcal{C}|}\left(g^{-1}\left(\frac{\mathbf{x}^*_{\mathcal{C}}}{r}\right)-\mathbf{\Sigma}_{\mathcal{C};\mathcal{U}}\mathbf{\Sigma}^{-1}_{\mathcal{U};\mathcal{U}}g^{-1}\left(\frac{\mathbf{x}^*_{\mathcal{U}}}{r}\right);\mathbf{\Sigma}_{\mathcal{C}|\mathcal{U}}\right)\\
        &\times \phi_{|\mathcal{U}|}\left(g^{-1}\left(\frac{\mathbf{x}^*_{\mathcal{U}}}{r}\right),\mathbf{\Sigma}_{\mathcal{U}}\right)\prod_{j\in \mathcal{U}}g^{-1'}\left(\frac{x^*_j}{r}\right) r^{-|{\mathcal{U}}|}f_R(r)dr.
    \end{split}
\end{equation} 

Although only one-dimensional integral appears in \eqref{censorLik}, the integrand \replace{includes }a $|\mathcal{C}|$-dimensional Gaussian distribution function. When approximating the integral using standard quadrature or Monte Carlo methods, one needs to compute $\Phi_{|\mathcal{C}|}$ for each sample point taken on $(1,r^*)$. This is only feasible when the number of locations $D$ is moderate. Additionally, this calculation will have to be repeated for each time replicate.

To avoid the integrating the process below the threshold, one could instead think of $X^*(\boldsymbol{s})$ as latent and draw from it using Monte Carlo methods. Consequently there is no need to compute the awkward likelihood \eqref{censorLik}. However, to update the Markov chain each time, it is now necessary to draw $\boldsymbol{x}^*_{\mathcal{C}}$ from a high-dimensional truncated distribution, which might again be computationally intensive.

Therefore, we propose to make a slight adjustment to the model in \eqref{scalemixtureModel}. Our new model is markedly more amenable to higher-dimensional inference, yet it keeps hold of the joint tail decay rates attained in the original model \citep[e.g.][]{huser2017modeling, huser2017bridging}. Equivalently, our new model has non-trivial asymptotically dependent and asymptotically independent submodels with the transition taking place in the interior of the parameter space in the case of our modified version of \eqref{huserWadsworth}.

\section{Model} \label{sec:model}
\subsection{Construction}\label{ourModelSec}
We alter the models in  Section \ref{hybrid} by adding an independent measurement error term to each component,
\begin{equation}\label{ourModel}
    X(\boldsymbol{s}_i):=X^*(\boldsymbol{s}_i)+\epsilon_i=R\cdot g(Z(\boldsymbol{s}_i))+\epsilon_i, 
\end{equation}
where $\epsilon_i\stackrel{iid}{\sim}N(0,\tau^2), \;i=1,\ldots, D$, and distribution of $R$ and the link function remain the same. That is, we add a simple nugget effect to the smooth process $X^*(\boldsymbol{s}_i)$. When drawing the latent processes below the threshold, we can condition on the smooth process and simply update the noisy one. Because these error terms are independent of each other, there is only a univariate integral involved in the full conditional likelihood. Also, when we update the smooth process $X^*(\boldsymbol{s})$ given the noisy process $X(\boldsymbol{s})$, no truncation or censoring is present and it is much easier to sample from the corresponding likelihood. Section \ref{sec:bayes} contains more details on the Markov Chain Monte Carlo (MCMC) updating scheme, where we show how this small alteration can hugely facilitate inference.

\subsection{Dependence Properties}\label{sec:dep}
We begin with the model \eqref{huserWadsworth} from \citet{huser2017modeling}, modified as in \eqref{ourModel}. Recall that $g(Z(\boldsymbol{s}))$ is a stationary process with standard Pareto margins possessing asymptotic independence; i.e., $P(g(Z(\boldsymbol{s}))>x)=x^{-1}$ and
\begin{equation}
      P(g(Z(\boldsymbol{s}_i))>x, g(Z(\boldsymbol{s}_k))>x)=\mathcal{L}_Z(x)\cdot x^{-\frac{1}{\eta_Z(h)}},\quad i\neq k,
\end{equation}
where $\mathcal{L}_Z(x)$ is slowly varying at infinity, and $\eta_Z(h)=(1+\rho(h))/2<1$ for the Gaussian correlation $\rho(h)<1$. 

\begin{figure}
    \centering
    \includegraphics[width=0.45\linewidth]{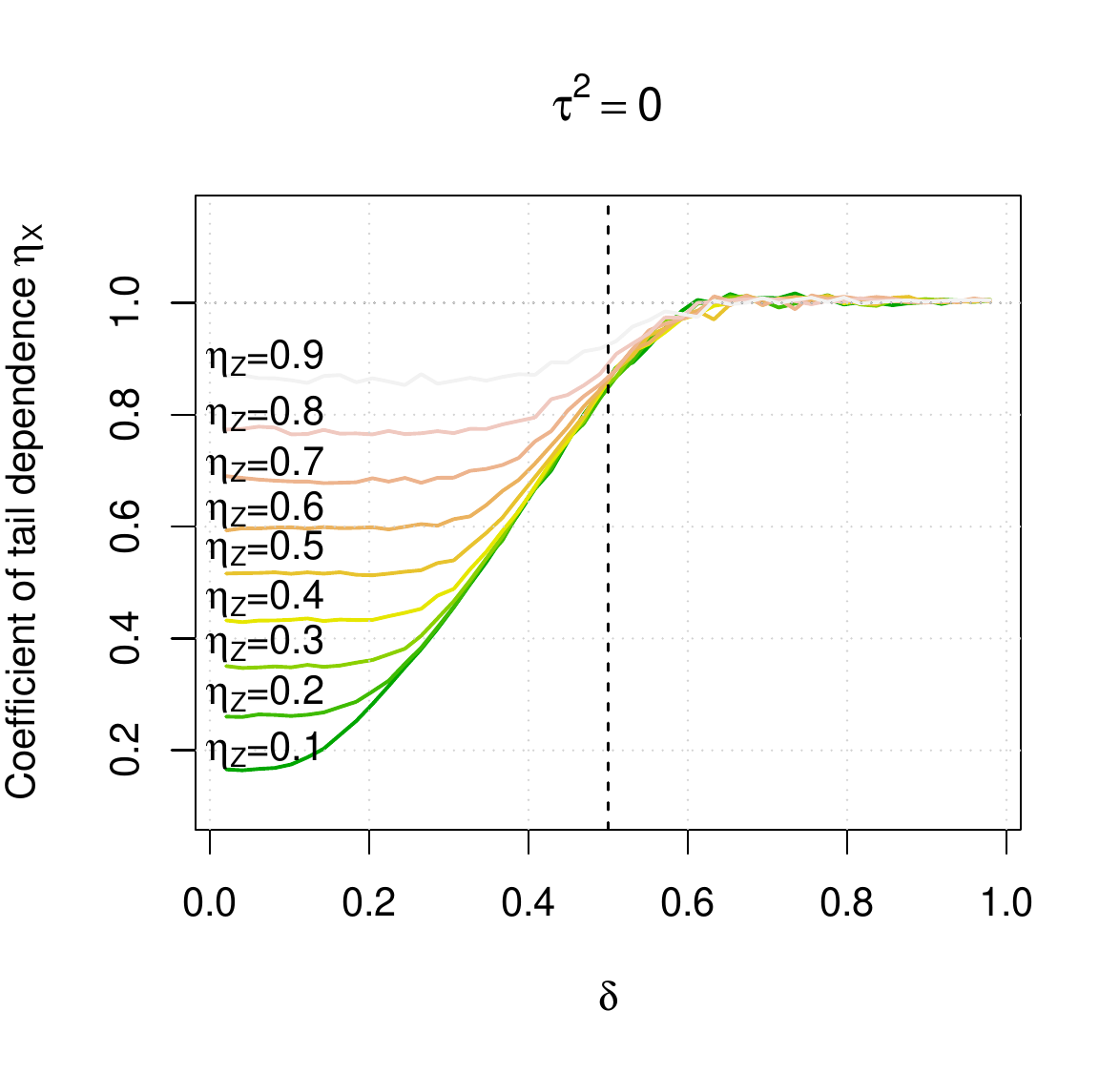}
    \includegraphics[width=0.45\linewidth]{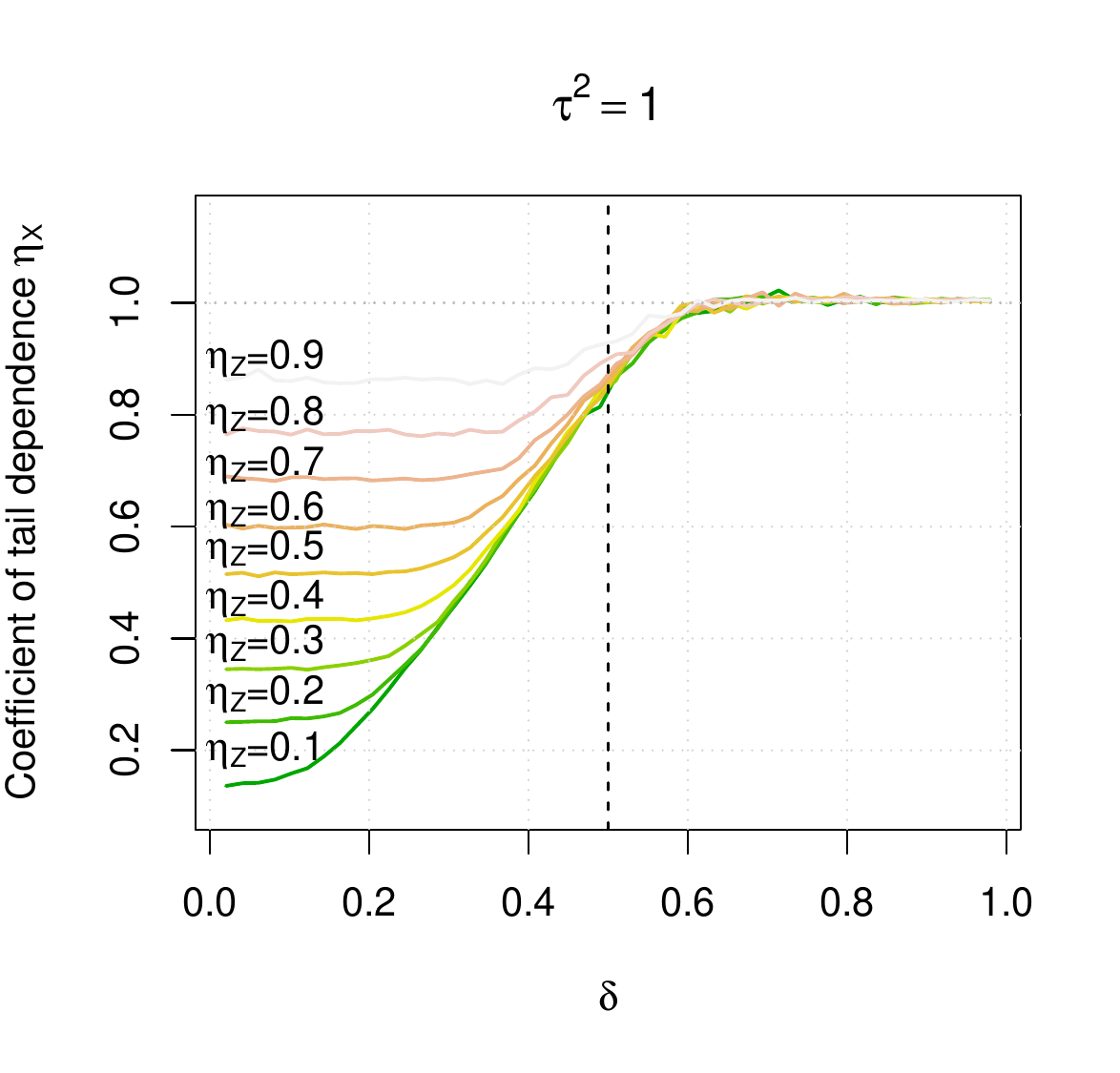}
    \caption{Coefficient of tail dependence approximated for the smooth Gaussian scale mixture processes ($\tau^2=0$) and the noisy processes ($\tau^2=1$) as a function of $\delta\in (0,1)$ for $\eta_W=0.1,\ldots,0.9$. The levels of dependence are similar for two models.}
    \label{numEta}
\end{figure}

Figure \ref{numEta} illustrates the estimated coefficient of tail dependence $\eta_X$ as a function of $\delta$ for $\eta_Z=0.1,\ldots,0.9$. For each combination of $\delta$ and $\eta_Z$, we generate 5,000,000 replicates from model \replace{\eqref{huserWadsworth}} (i.e. $\tau^2=0$) and model \eqref{ourModel} with $\tau^2=1$ respectively. For each replicate, we sample $(Z_i, Z_k)$ from a Gaussian copula with correlation $2\eta_Z-1$. We then numerically approximate the joint survival probability in \eqref{tailDependence} to obtain an estimate of $\eta_X$. The left panel of Figure \ref{numEta} clearly shows that the smooth transition from asymptotic independence to asymptotic dependence takes place around $\delta=1/2$, confirming the results from \citet{huser2017modeling} with a reasonable bias; the right panel shows that adding a measurement error has little effect on the tail dependence because $\eta_X$ exhibits similar behavior. This result invites investigation of whether the flexible asymptotic properties in \citet{huser2017modeling} are preserved in the altered model. In the following, we generalize the problem from the specific model of \citet{huser2017modeling} for the process $X^*$ in~\eqref{ourModel}, to any $X^*$ with a wide class of marginal tail behaviors.

The impact of the additive Gaussian nugget effect on the extremal dependence of $X$ depends upon the marginal tail heaviness of $X^*$: roughly, the heavier the tail of $X^*$, the less the impact of the noise. Since we take a copula-like approach and employ $X$ as the spatial dependence model, we assume its margins, and those of $X^*$, are identical over space.

We will focus on the dependence of $X$ under two broad classes of marginal distribution for $X^*$: regularly varying tails, and Weibull-like tails. A measurable function $g:\mathbb{R}_+ \to \mathbb{R}_+$ is said to be regularly varying at infinity with index $\kappa$ if for all $x>0$
$\lim_{t\to\infty}g(tx)/g(t)=x^\kappa$; we write $g \in \RV_\kappa$. When $\kappa=0$, the function is slowly varying.

Regularly varying tails are defined through the survival function being regularly varying at infinity, i.e., $\P(X^*>x) \in \RV_{-\alpha}$, $\alpha>0$. Weibull-like tails are defined through the survival function
\begin{align}
 \P(X^* > x) \sim u(x)\exp(-\theta x^\alpha ), \qquad \theta, \alpha>0, \qquad x \to \infty \label{eq:lightsf}
\end{align}
where $u \in \RV_{\kappa}$, and $\alpha$ is termed the Weibull index. We also assume that $X^*$ has a density satisfying $f_{X^*}(x) \sim v(x)\exp(-\theta x^\alpha)$, with $v(x) = u(x)(\theta\alpha x^{\alpha-1})$. The main results are now summarized in Proposition~\ref{prop:dep}. 

\begin{proposition}
\label{prop:dep}
With definitions and notation as above:
\begin{enumerate}
 \item If $X^*$ has a regularly varying tail, or Weibull-like with Weibull index $\alpha<1$, then $\chi_{X}=\chi_{X^*}$ and $\eta_X = \eta_{X^*}$. \label{item1}
 \item If $X^*$ has a Weibull-like tail with Weibull index $\alpha=1$ then 
 \[\chi_{X}\in [\E(e^{\theta\min(\epsilon_1,\epsilon_2)})/\E(e^{\theta \epsilon}),\E(e^{\theta\max(\epsilon_1,\epsilon_2)})/\E(e^{\theta \epsilon})]\chi_{X^*}\] \label{item2}
 and $\eta_X = \eta_{X^*}$. Note if $\chi_{X^*}=0$ then so is $\chi_{X}$.
 \item If $X^*$ has a Weibull-like tail with Weibull index $\alpha>1$ then: \label{item3}
 \begin{enumerate}
  \item If $\alpha \in(1,2)$, then $\eta_X = \eta_{X^*}$ \label{item3a}
  \item If $\alpha=2$, then an interval can be given for $\eta_X$ (see Expression~\eqref{eqn:a2etarange}).
  \item If $\alpha>2$, then $\eta_X = 1/2$.
 \end{enumerate}
 \end{enumerate}
\end{proposition}
The proof of Proposition~\ref{prop:dep} is given in Appendix~\ref{app:proofs}. 

For the process of \citet{huser2017modeling}, described in equation~\eqref{huserWadsworth}. In this case, $X^*$ always has a regularly varying tail, so Proposition~\ref{prop:dep} part~\ref{item1} gives $\eta_{X}(h)=\eta_{X^*}(h)$, and $\chi_{X}(h)=\chi_{X^*}(h)$, with $\eta_{X^*}(h)$ and $\chi_{X^*}(h)$ given in Section \ref{hybrid}. This means the flexible asymptotic properties in \citet{huser2017modeling} are preserved in the altered model.

Another popular model for spatial data is the $t$-process \citep{roislien-2006a}, which is a Gaussian scale mixture for which the mixing variable $R^2$ follows an inverse gamma distribution. If $X^*$ is a $t$-process, then it is asymptotically dependent with a regularly varying tail, so $\eta_{X}(h)=\eta_{X^*}(h) = 1$ and $\chi_{X}(h)=\chi_{X^*}(h)>0$ by Proposition~\ref{prop:dep}. Similarly, the skew-$t$ process \citep{padoan-2011a,morris-2017a} is regularly varying and asymptotically dependent, so the same conclusions apply.

For the Gaussian scale mixture of \citet{huser2017bridging}, $X^*$ either has regularly varying or Weibull-like tails depending on the distribution of the scaling variable $R$. In particular when $\beta>0$ in distribution~\eqref{hot}, $X^*$ has a Weibull-like tail with Weibull index $\alpha=2\beta/(\beta+2) <2$. As such, parts~\ref{item1},~\ref{item2} and~\ref{item3a} of Proposition~\ref{prop:dep} are relevant. In all cases $\eta_{X^*}(h) = \eta_X(h) = \{(1+\rho(h))/2\}^{\beta/(\beta+2)}$, and $\chi_{X}(h)=\chi_{X^*}(h)=0$ for $\beta>0$. When $\beta=0$ then $X^*$ has a regularly varying tail, and is asymptotically dependent with $\eta_{X}(h)=\eta_{X^*}(h)=1$, and 
\begin{align*}
\chi_{X^*}(h)=\chi_{X}(h) = 2\{1-T_{\gamma+1}((\gamma+1)^{1/2}(1-\rho(h))/(1-\rho(h)^2)^{1/2})\}>0,
\end{align*}
where $T_{\nu}$ is the cdf of the Student-$t$ distribution with $\nu$ degrees of freedom.

We note that the process $X^*$ in equation~\eqref{scalemixtureModel} is constructed only for its dependence properties, and there is no ``natural" scale on which to express it. For example, considering the process of \citet{huser2017modeling} we could also write
\begin{align}
X^*(s) = E + V(s), \label{eqn:hwE}
\end{align}
with $E \given \delta \sim \mbox{Exp}(\delta/(1-\delta))$, $V(s)=\log g(Z(s))$, which has the same dependence structure as defined in~\eqref{scalemixtureModel} and~\eqref{huserWadsworth}, since it is obtained through a monotonic marginal transformation. Taking $X^*$ from~\eqref{eqn:hwE}, Proposition~\ref{prop:dep} part~\ref{item2} gives $\eta_{X}(h)=\eta_{X^*}(h)$. Asymptotic dependence of $X^*$ implies asymptotic dependence of $X$, but only bounds on $\chi_{X}(h)$ are available. 

In practice, if choosing a marginal scale for $X^*$, there may be a trade-off between theoretical desires and computational practicality. Supposing that we wish $X$ to inherit the properties of $X^*$, a heavy-tailed choice is best.

\section{Bayesian Inference}\label{sec:bayes}
\subsection{Hierarchical Model}
We define a Bayesian hierarchical model based on the process \eqref{ourModel} defined in Section \ref{ourModelSec} and use a MCMC algorithm to fit to the data. For the reasons outlined in Section \ref{censored}, we assume data are censored below a high threshold $u$. In \eqref{ourModel}, the mixing parameter \replace{$\btheta_R$} controls both joint and marginal behavior of the response $X$, which we would prefer to separate. Therefore, motivated by the theory of univariate extremes, we first assume our observations above the same high threshold $u$ are generalized Pareto distributed, and we include a marginal transformation in the hierarchical model. 

\edit{Let $\{Y(\boldsymbol{s}): \boldsymbol{s}\in \mathcal{S}\}$ denote  the observed process}. We define a marginal transformation $T(y)$ as follows:
\begin{equation}\label{marginal_Transform}
    T(Y(\boldsymbol{s}))=F_{X|\btheta_R,\tau^2}^{-1}\circ F_{Y|u,\sigma,\xi}(Y(\boldsymbol{s})),
\end{equation}
where $F_{X|\btheta_R,\tau^2}$ is the marginal distribution function for process \eqref{ourModel}, and 
\begin{equation}
    F_{Y|u,\sigma,\xi}(y)=\left\{
                \begin{array}{ll}
                 p,& y\leq u,\\
                 p+(1-p)F_{GPD|u,\sigma,\xi}(y),&y>u,
                \end{array}
              \right.  
\end{equation}
where $p=P(Y(\boldsymbol{s}_i)\leq u)$, $u$ is a high threshold, and $F_{GPD|u,\sigma,\xi}(y)=1-[1+\xi(y-u)/\sigma]^{-1/\xi}$, \edit{with support $\{y\geq u:1+\xi(y-u)/\sigma\geq0\}$}. Conditioning on the smooth process $X^*(\boldsymbol{s}_i)=x^*(\boldsymbol{s}_i)$, which was not truncated, the censored likelihood for an observation $y(\boldsymbol{s}_i)$ can be derived as
\begin{equation}\label{fullY}
    \varphi(y(\boldsymbol{s}_i) \given \boldsymbol{X}^*,\tau^2,\btheta_R,\boldsymbol{\theta}_{GPD},p)= \left\{
                \begin{array}{ll}
                  \Phi\left(\frac{F^{-1}_{X|\btheta_R,\tau^2}(p)-x^*(\boldsymbol{s}_i)}{\tau}\right) &\text{if } y(\boldsymbol{s}_i)\leq u,\\
                \frac{1}{\tau}\phi\left(\frac{T(y(\boldsymbol{s}_i))-x^*(\boldsymbol{s}_i)}{\tau}\right)\cdot\frac{f_{Y|u,\sigma,\xi}(y(\boldsymbol{s}_i))}{f_{X|\btheta_R,\tau^2}(T(y(\boldsymbol{s}_i)))}&\text{if } y(\boldsymbol{s}_i)>u,
                \end{array}
              \right.
\end{equation}
where $\boldsymbol{\theta}_{GPD}=(\sigma,\xi)$. Note that there are only univariate calculations required in \eqref{fullY}, compared to \eqref{censorLik} for which we have to estimate the $|\mathcal{C}|$-dimensional Gaussian distribution functions. In addition, since $Y_i$ and $Y_k$ are independent conditioning on the smooth process ($i\neq k$), the joint likelihood of the whole vector $\boldsymbol{y}$ is simply $f_{\boldsymbol{Y}}(\boldsymbol{y})=\prod_{i=1}^D \varphi(y(\boldsymbol{s}_i) \given \boldsymbol{X}^*,\tau^2,\boldsymbol{\theta}_{R},\boldsymbol{\theta}_{GPD},p)$. Likelihoods for independent time replicates are simply multiplied together, and the proportion of censored observations $p$ can be treated as a known parameter or an unknown parameter that enters the hierarchical model. See Appendix \ref{hierarch_mod} for a complete statement of the hierarchical model. The priors for the model parameters are
\begin{align}\label{hierarchical}
        \tau^2&\sim \text{IG}(\alpha,\beta), &\sigma & \sim \text{halfCauchy}(1),& \xi & \sim \text{Unif}(-0.5,0.5),
\end{align}
where halfCauchy(1) refers to the positively truncated standard Cauchy distribution. \edit{We implement this methodology for two different Gaussian scale mixture models: that of \citet{huser2017modeling}, and that of \citet{huser2017bridging}}. The priors for $\boldsymbol{\theta}_R$ are $\delta\sim U(0,1)$ for the former, and $\beta\sim \text{halfCauchy}(1)$ for the latter. The prior for $\btheta_C$ depends on the choice of the covariance function. In our implementations, we adopt the Mat\'{e}rn covariance function with $\btheta_C=(\rho,\nu)$, where $\rho$ is the range parameter, and $\nu$ is the smoothness parameter. The prior is then set to subject to two independent half Cauchy distributions with scale parameter 1.

\subsection{Gibbs Sampler} \label{gibbs}
To estimate the posterior distribution of the model parameters, we apply random walk Metropolis (RWM) algorithm using Log-Adaptive Proposals (LAP) as our adaptive tuning strategy \citep{shaby2010exploring}. Since conjugate priors are not available, we use random walk Metropolis-Hastings update steps.

At each MCMC iteration, we first update the smooth process $X^*$ conditioning on the true values for all non-censored sites, current values for $X^*$ and all other model parameters:
\begin{equation}\label{xstar}
    \varphi(\boldsymbol{X}^*\given \cdots) \;\propto\; \prod_{i=1}^D\varphi(y_i\given \boldsymbol{X}^*, \tau^2,\btheta_R,\boldsymbol{\theta}_{GPD},p)\cdot\varphi(\boldsymbol{X}^*|R,\btheta_C),
\end{equation}
where the likelihood function of $X^*$ conditioning on the random scaling factor $R$ is calculated in Appendix \ref{hierarch_mod}. We then update $R$ using its conditional posterior distribution
\begin{equation*}
    \varphi(R\given \cdots)\;\propto\;\varphi(\boldsymbol{X}^*\given R,\btheta_C)\cdot\varphi(R\given \btheta_R).
\end{equation*}

Since time independence is assumed, we can update $\boldsymbol{X}^*_t$ and $R_t$ in a parallel fashion across $t=1, \ldots, T$. The other parameters are updated similarly using adaptively-tuned random walk Metropolis-Hastings updates, with the likelihood~\eqref{fullY} multiplied by the corresponding priors in \eqref{hierarchical}.

\section{Simulation Studies}\label{sec:simu}
In this section, we present simulation results and conduct coverage analysis to investigate, firstly, whether the MCMC procedure is able to draw accurate inference on model parameters, and secondly, in the case of the modified version of model \eqref{huserWadsworth}, to check whether our model captures asymptotic dependence characteristics correctly even when the data-generating model is different from the fitted model.

\subsection{Parameter Estimation}\label{paramEstimate}
To verify the accuracy of inference made by MCMC sampling, we generate data from model \eqref{ourModel} in Section \ref{ourModelSec}, in both the special case of the \citet{huser2017bridging} model \eqref{hot} and the special case of the \citet{huser2017modeling} model \eqref{huserWadsworth}, with the addition of nugget terms.  In both cases, we use $D=200$ sites uniformly drawn from the unit square $[0,1]^2$, with the latent Gaussian processes $Z(\boldsymbol{s})$ are generated using a Mat\'ern covariance with smoothness parameter $\nu=3/2$.  The characteristic length scale parameter is set to $\rho=0.05$ in the case of model \eqref{hot} and $\rho=0.1$ in the case of model \eqref{huserWadsworth}.  

For model \eqref{hot}, we use $T=20$ independent temporal replications, and set $\beta=0.5$, $\gamma = 1$ (which is fixed during estimation), and nugget variance $\tau^2 = 0.2^2$.  This represents a challenging case, with a fairly long tail and nugget that is small relative to the scale of $X^*(\boldsymbol{s})$. For model \eqref{huserWadsworth}, we use $T=40$ independent temporal replications, and set the nugget variance to $\tau^2 = 3^2$.  Because the latent $Z(\boldsymbol{s})$ is transformed to Pareto, $\tau^2 = 3^2$ is still small compared to the scale of smooth process $X^*(\boldsymbol{s})$; see the Supplementary Material for a more in-depth discussion on the effects of $\tau^2$. We consider two different scenarios for the dependence parameter $\delta$: $\delta=0.3$ and $\delta=0.7$, corresponding to asymptotic independence and asymptotic dependence respectively.  Finally, in all cases the processes are marginally transformed to generalized Pareto distribution with $(u,\sigma,\xi)=(11,1,0)$, where $u$ and $p=0.8$, the proportion of censored observations, are treated as known parameters.

The attenuation constants used in the LAP algorithm are $c_0=10, c_1=0.8$. The prior for $\rho$ is halfCauchy$(1)$, and the priors for the other parameters are specified in \eqref{hierarchical}, where $(\alpha,\beta)=(0.1,0.1)$ so that the prior for $\tau^2$ is fairly noninformative.  We ran each MCMC chain for 400,000 iterations and thinned the results by a factor of 10. The parallelism of updating $\boldsymbol{X}^*_t$ and $R_t$ is implemented in \texttt{R} via the \texttt{foreach} routine with \texttt{doParallel} package as a backend \citep{analytics1weston}.

\subsection{Coverage Analysis}
We now study the coverage properties of the posterior inference based on the MCMC sampler for the posterior credible intervals with 100 simulated datasets drawn from the \citet{huser2017bridging} and \citet{huser2017modeling} models, under each of the scenarios described in the previous section. 

\begin{figure}
 \centering
  \begin{minipage}{0.5\textwidth}
    \centering
    $\beta = 0.5$ \\ \vspace{1em}
    \includegraphics[width=0.49\linewidth]{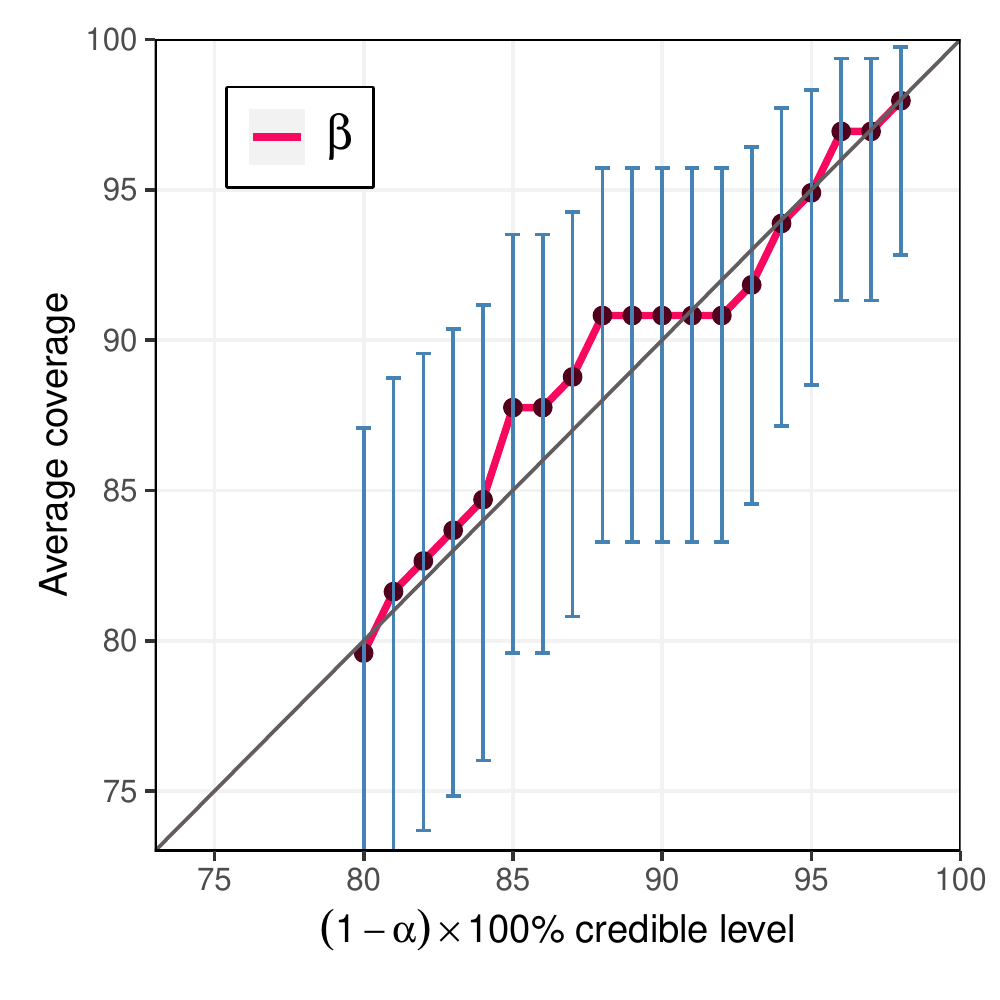}
    \includegraphics[width=0.49\linewidth]{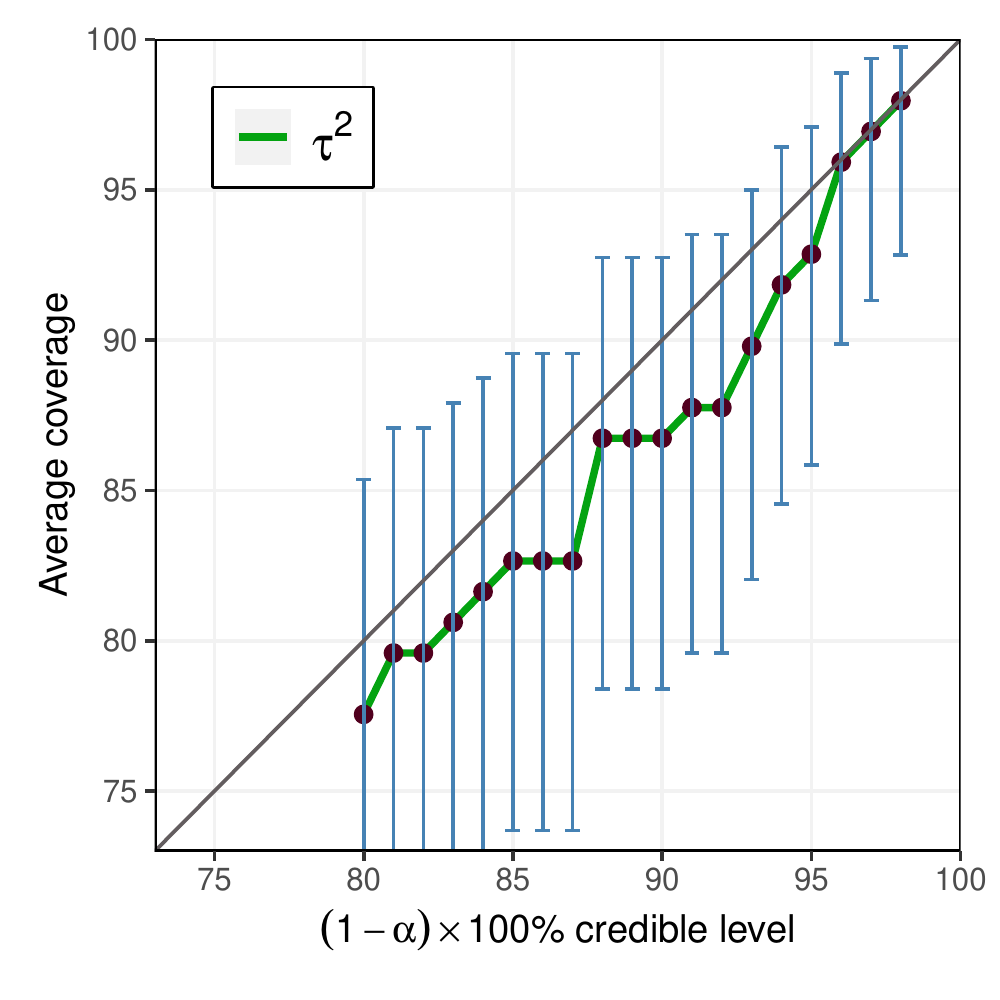}
  \end{minipage}
    \caption{Empirical coverage rates of credible intervals for the \citet{huser2017bridging} model: $\beta=0.5,\; \tau^2=0.2^2$. The error bars are 95\% binomial confidence intervals for the coverage probability.}
    \label{Coverage_hot}
\end{figure}

\begin{figure}
  \begin{minipage}{0.5\textwidth}
    \centering
    $\delta = 0.3$ \\ \vspace{1em}
    \includegraphics[width=0.49\linewidth]{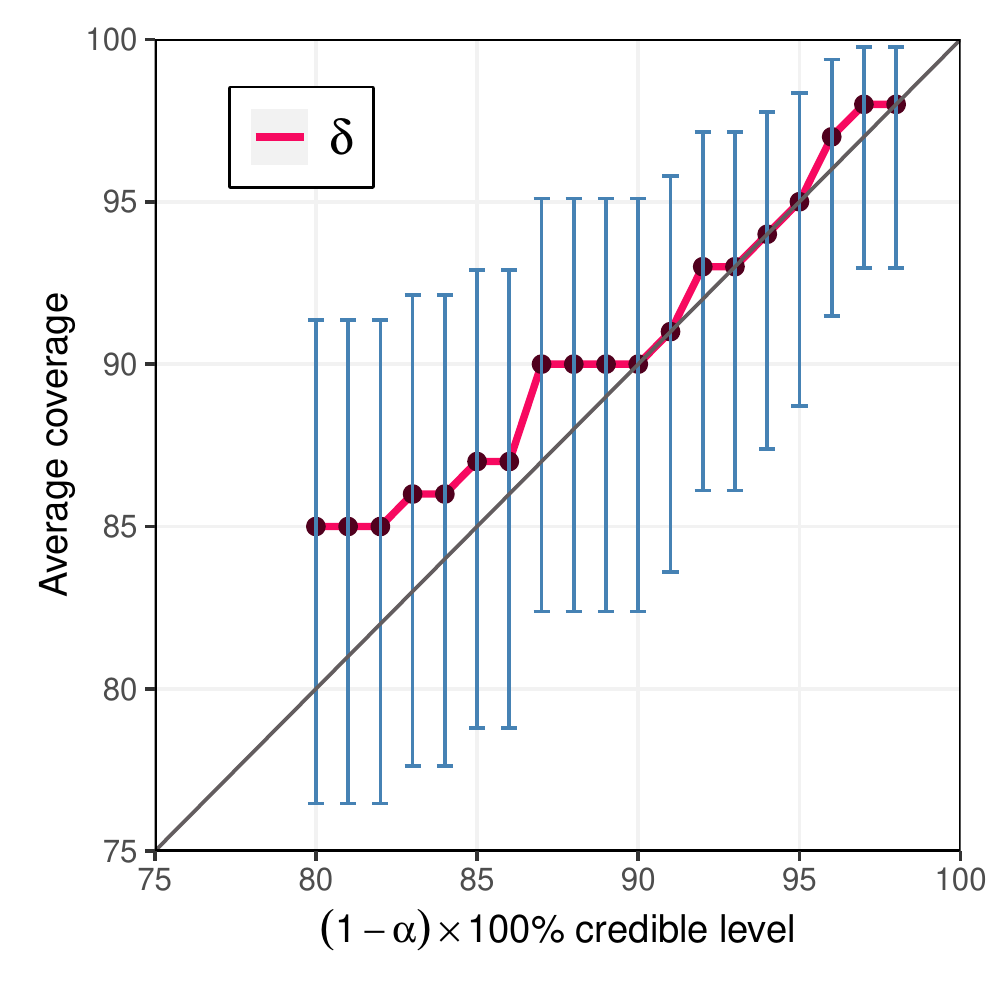}
    \includegraphics[width=0.49\linewidth]{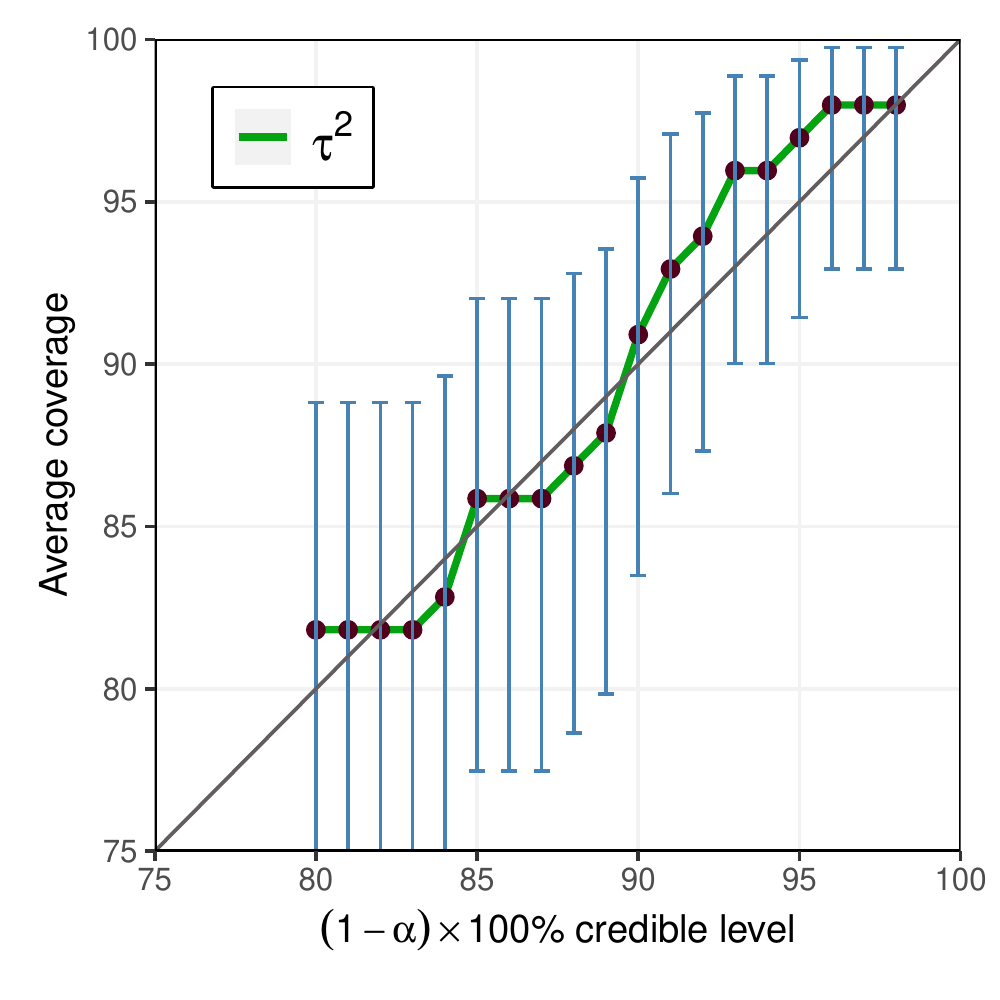}
  \end{minipage}
  \begin{minipage}{0.5\textwidth}
    \centering
    $\delta = 0.7$ \\  \vspace{1em}
    \includegraphics[width=0.49\linewidth]{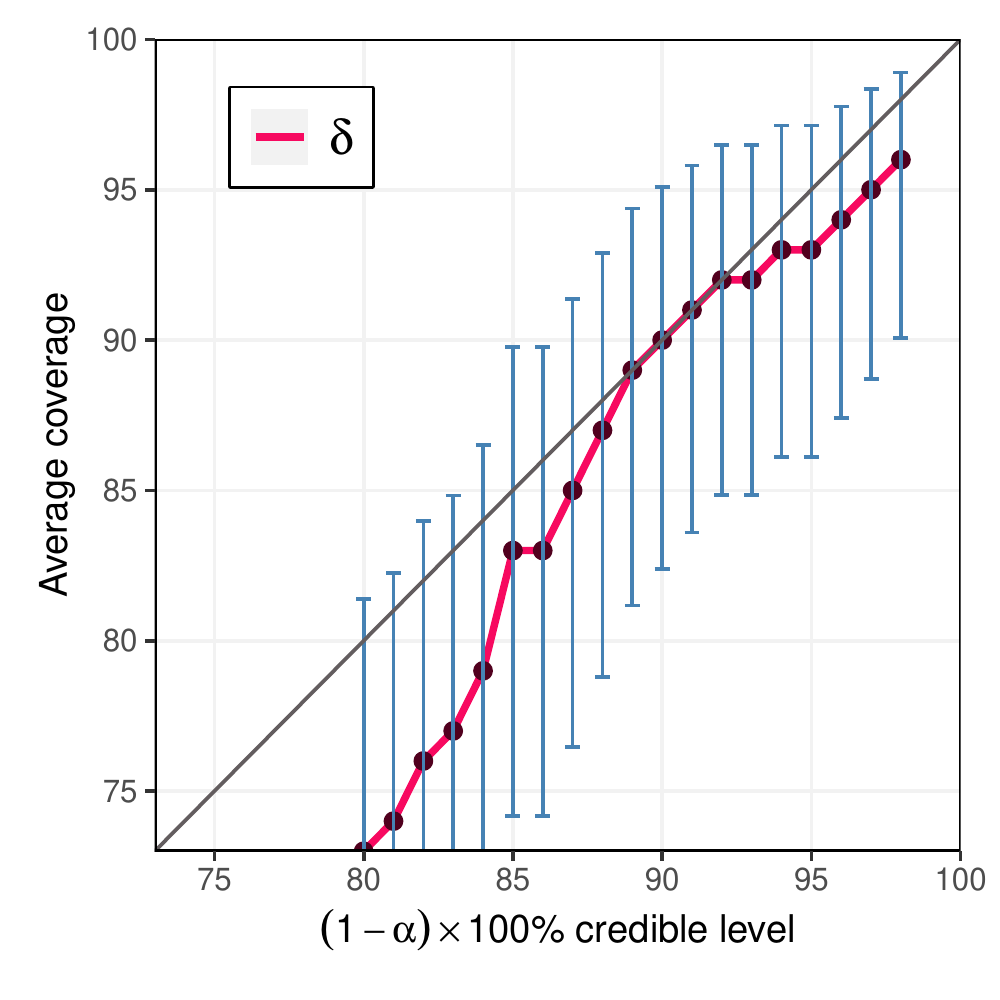}
    \includegraphics[width=0.49\linewidth]{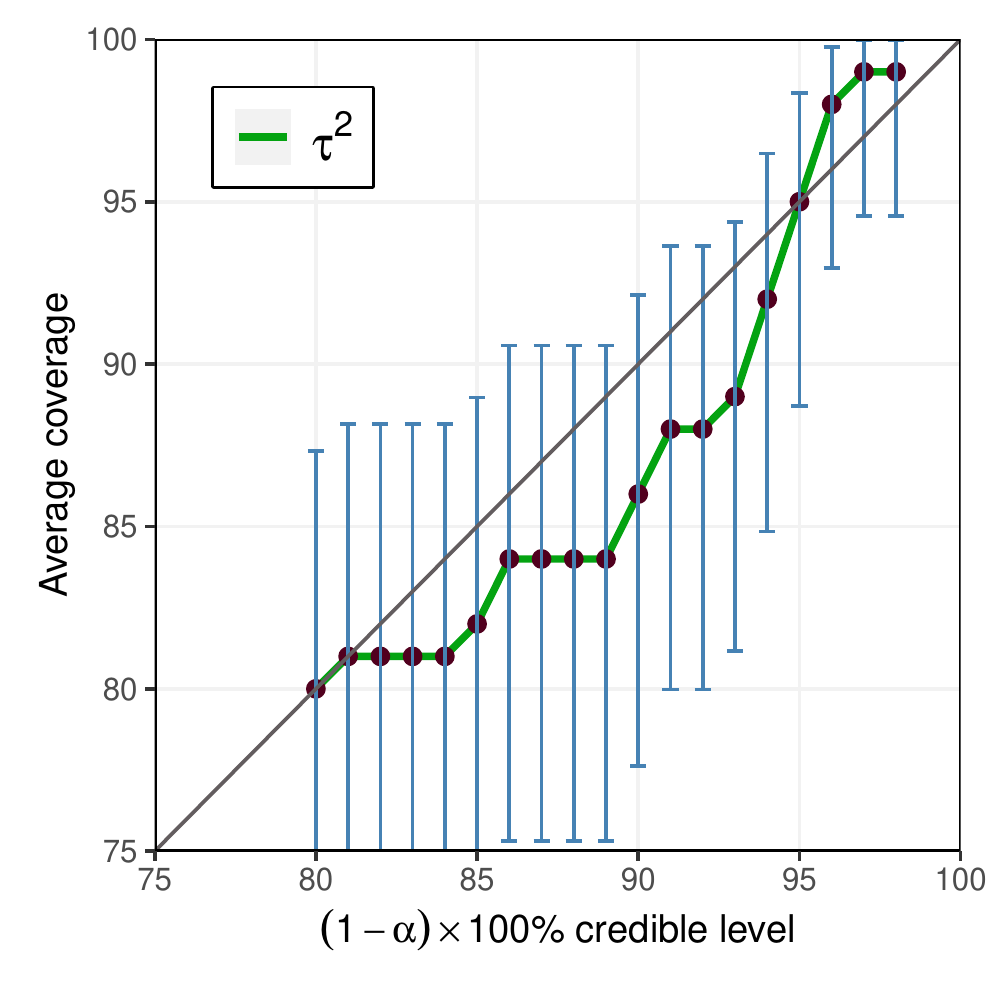}
  \end{minipage}
    \caption{Empirical coverage rates of credible intervals for the \citet{huser2017modeling} model under two designs: $\delta=0.3$ (left) and $\delta=0.7$ (right). For both designs, $\tau^2=3^2$. The error bars are 95\% binomial confidence intervals for the coverage probability.}
    \label{Coverage}
\end{figure}

Figures \ref{Coverage_hot} and \ref{Coverage} shows the empirical coverage rates of highest posterior density credible intervals of several sizes, along with standard binomial confidence intervals. In all cases, we can see that the sampler performed well in generating posterior inference that is well calibrated, with close to nominal frequentist coverage. The coverage for larger $\delta$ is may be slightly different than nominal for large $\alpha$, but overall the results are quite good.

\subsection{Simulation with Mis-specified Models}\label{mis-specified}
We now fit our model to data generated from other distributions to validate its ability to capture the tail dependence characteristics under mis-specification.  We simulate datasets from models referenced in Section \ref{hybrid}, and use the sampler described in Section \ref{gibbs} to fit model \eqref{ourModel} with the \citet{huser2017modeling} latent process. The data were generated using four different simulation designs:
\begin{itemize}
    \item Skew-$t$ process from \citet{morris2017space} with $(a,b)=(6,16)$ (asymptotically dependent);
    \item Gaussian scale mixture from \citet{huser2017bridging} with $\beta=0$ (asymptotically dependent);
    \item Gaussian scale mixture from \citet{huser2017bridging} with $\beta=1\text{ or }5$ (asymptotically independent).
\end{itemize}
For each simulation design, we simulate a single dataset using $D=100$ locations uniformly distributed on $[0,1]^2$ and $T=40$ independent time replicates. The Mat\'ern covariance function with $\nu=3/2$ and $\rho=1$ is again specified for the latent Gaussian processes. For the skew-$t$ process, $R_t\sim \text{IG}(3,8)$ to give a $t$ distribution with $6$ degrees of freedom, and $\lambda=3$ to simulate moderate skewness. For the last two designs, the $R_t$ were generated as described in \eqref{hot}, with $\gamma=1$.

To obtain good starting values for the latent smooth process $\boldsymbol{X}^*$ for MCMC, we first marginally transform the simulated data to noisy scale mixture variables $\boldsymbol{X}$ independently at each location using the following procedure. Following the semi-parametric procedure of \citet{coles1991modelling}, we estimate each marginal distribution as a blend of the generalized Pareto distribution function above a high marginal threshold, and the empirical distribution function below that threshold. Fixing initial values for $(\delta,\tau^2)$, we then transform the margins to noisy scale mixtures via $X_{jt}=F_{X|\delta,\tau^2}^{-1}\circ\hat{F}_{\boldsymbol{s}_j}(Y_{jt})$. The next step is to run a Metropolis algorithm using the full conditional distribution $\varphi(\boldsymbol{X}^*\given \cdots)$ (see \eqref{xstar}) 100 times and save the last random walk states as initial values for $\boldsymbol{X}^*$. This procedure is also used for data analysis in Section \ref{data_analysis}. Finally, with initial values in hand, the datasets from each design were fit using a fully Bayesian approach that simultaneously updates marginal and spatial dependence parameters. 

\begin{figure}
    \centering
    \begin{subfigure}[b]{0.4\linewidth}
    \centering
    \caption{\textbf{Skewed-$t$ process}}
    \includegraphics[width=\linewidth]{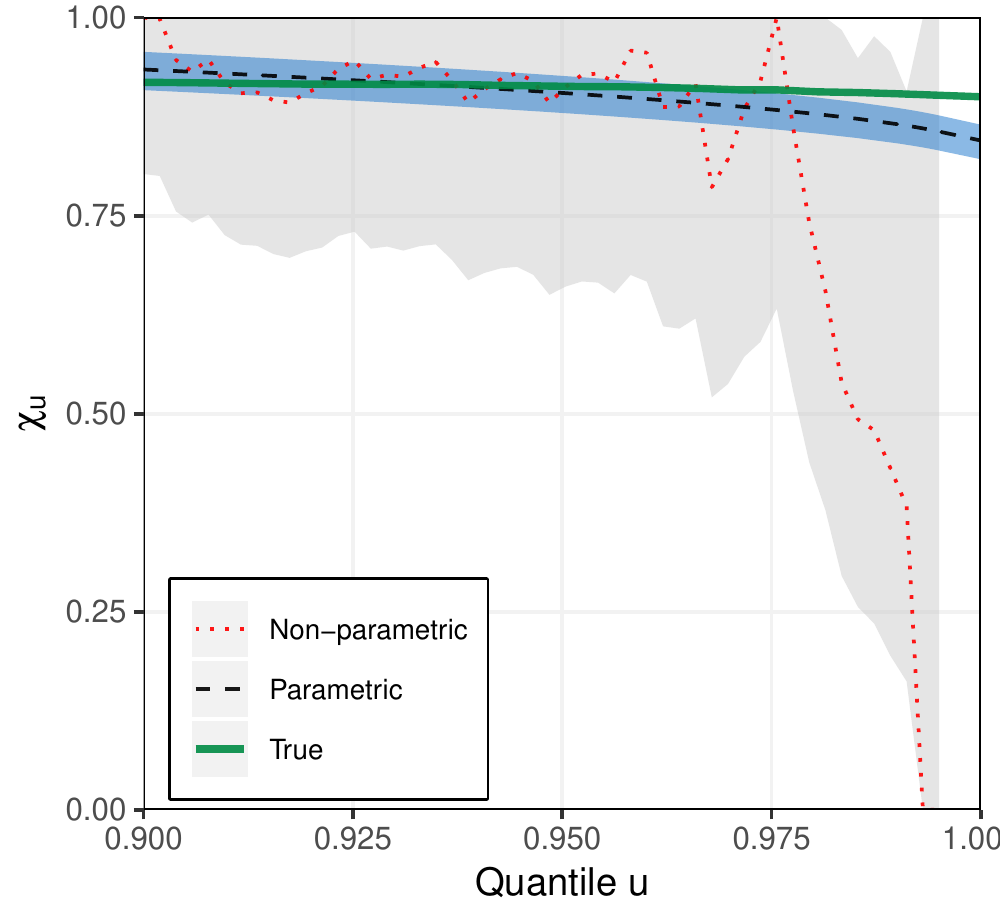}
    \end{subfigure}
    \begin{subfigure}[b]{0.4\linewidth}
    \centering
    \caption{\textbf{\citet{huser2017bridging} with $\beta=0$}}
    \includegraphics[width=\linewidth]{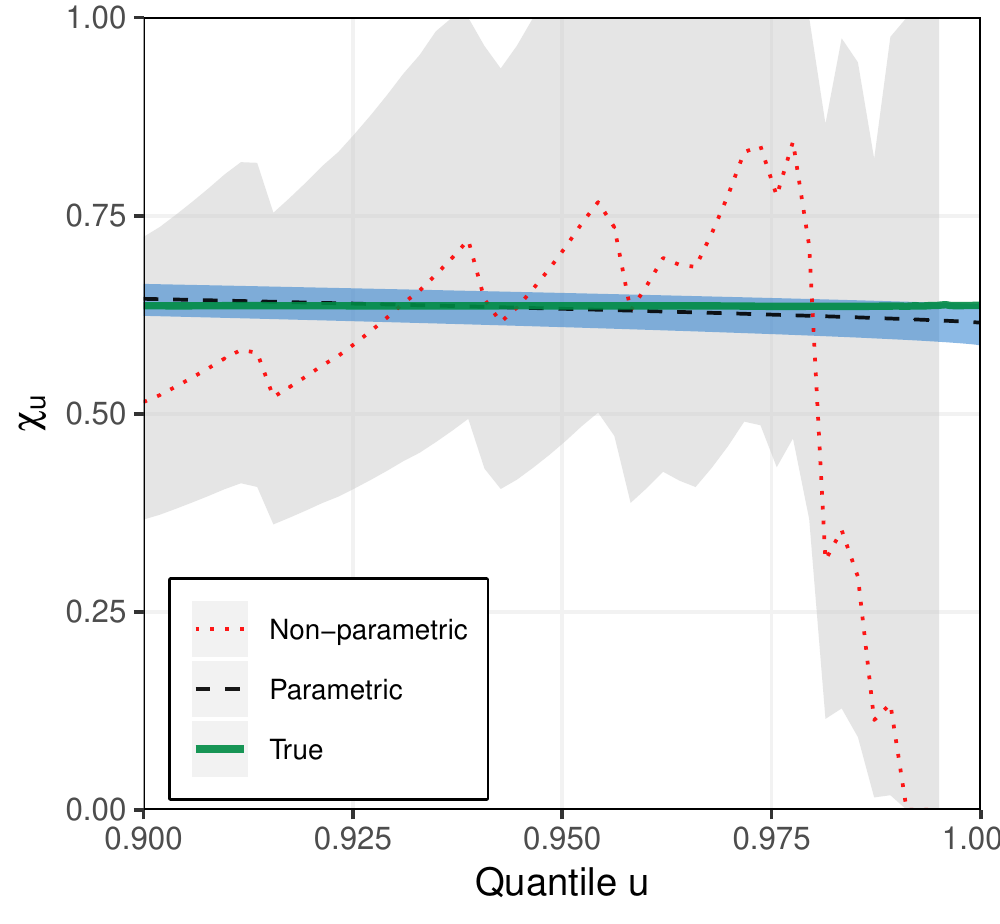}
    \end{subfigure}
    
    \begin{subfigure}[b]{0.4\linewidth}
    \centering
    \caption{\textbf{\citet{huser2017bridging} with $\beta=1$}}
    \includegraphics[width=\linewidth]{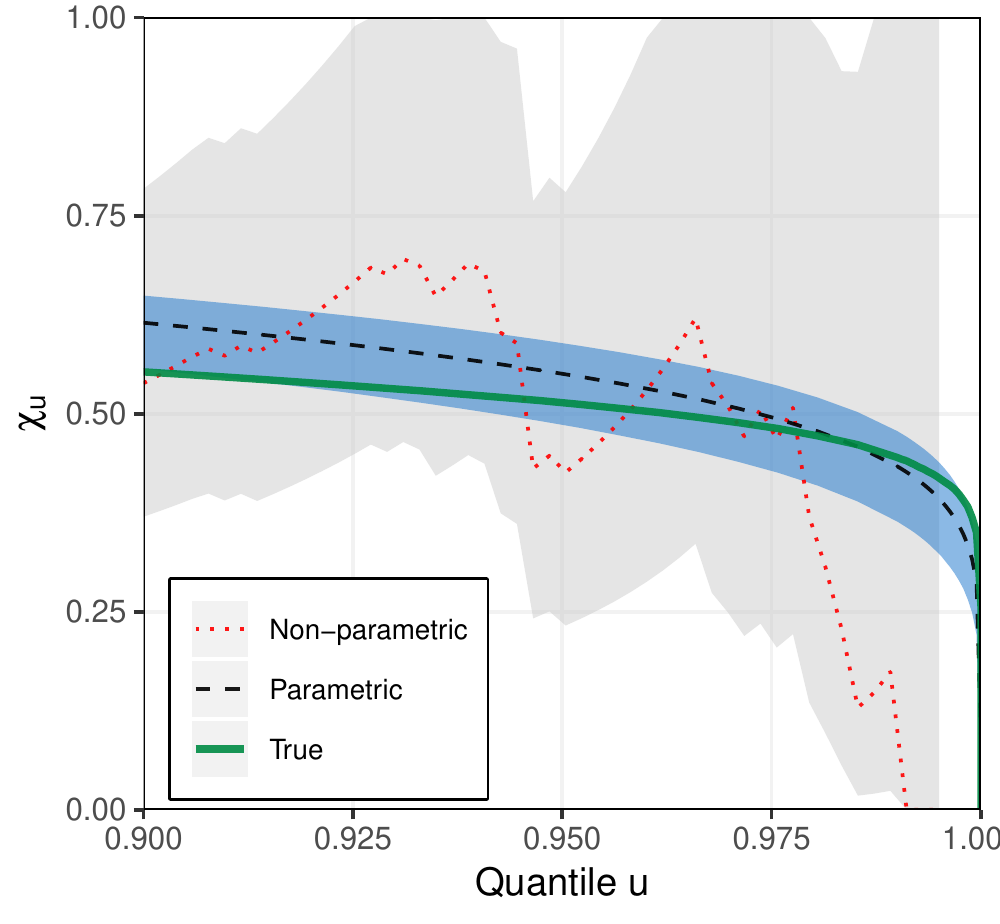}
    \end{subfigure}
    \begin{subfigure}[b]{0.4\linewidth}
    \centering
    \caption{\textbf{\citet{huser2017bridging} with $\beta=5$}}
    \includegraphics[width=\linewidth]{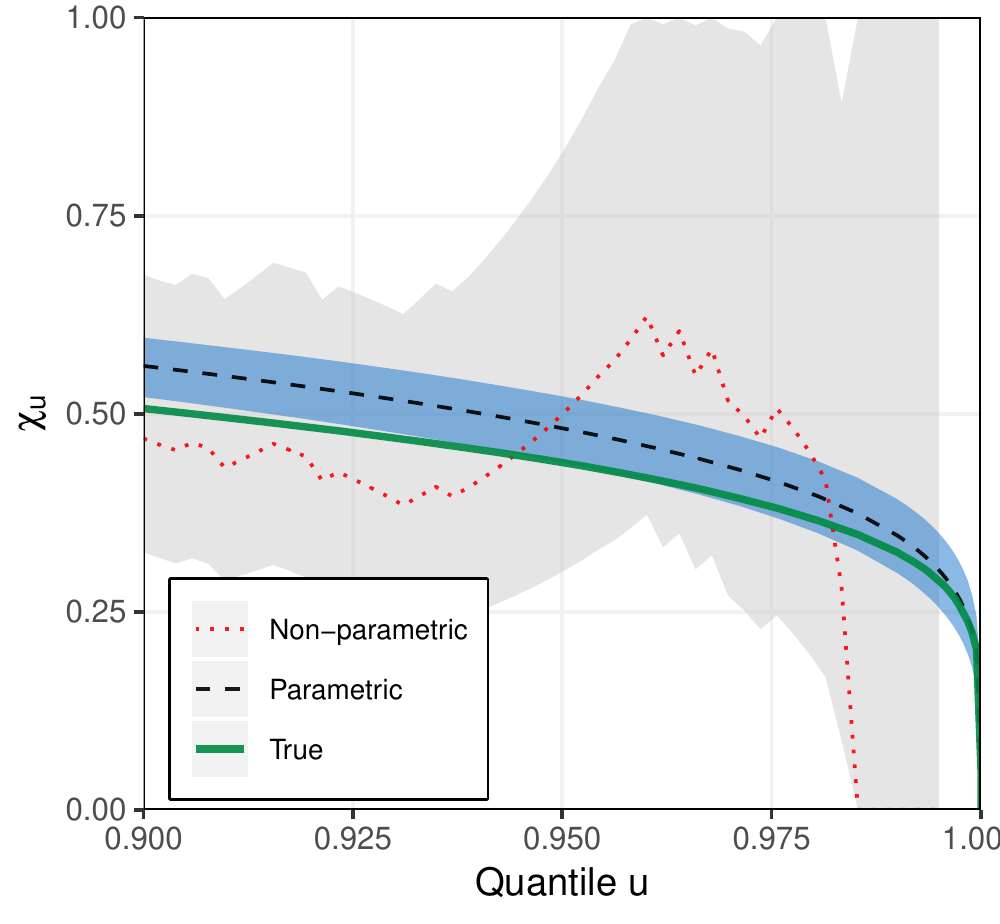}
    \end{subfigure}
    \caption{Estimated coefficients $\chi_u(h)$, $u\in[0.9,1]$, for two points at distance $||\boldsymbol{s}_1-\boldsymbol{s}_2||=1$, using data simulated from the mis-specified models with Mat\'ern correlation function with $\nu=3/2$ and $\rho=1$. The two scenarios in the top row are asymptotically dependent, while the scenarios in the bottom row are asymptotically independent. Each simulated data set has 100 uniform locations in $[0,1]^2$, and 40 time replicates. Solid green lines show true $\chi_u(h)$ function. The black dashed lines show the averaged curve from the posterior samples of fitted model, while the blue shaded areas are 95\% credible envelopes. The red dashed lines show the nonparametric estimates from the simulated data sets, while the grey shaded areas are pointwise 95\% binomial confidence intervals.}
    \label{chi_enve}
\end{figure}

Figure \ref{chi_enve} displays the nonparametric and model-based estimates of the upper tail dependence $\chi_u(h)$ defined in \eqref{chiu}. To generate nonparametric estimates of $\chi_u(h)$ at distance $h=\|\boldsymbol{s}_1-\boldsymbol{s}_2\|$, we look at all pairs of points whose locations are $h$ apart (within some small $\varepsilon$ tolerance), and compute the ratio of empirical probabilities $\hat{\chi}_u(h)$. This is similar to an empirical variogram estimator. The nonparametric confidence envelopes are obtained by computing pointwise binomial confidence intervals, pretending that each pair of points is independent from each other pair of points. For parametric estimates, we take samples from the converged MCMC chain, and use parameters from each iteration to simulate $10^8$ pairs of $(X(\boldsymbol{s}_1), X(\boldsymbol{s}_2))$ based on our model to generate a smooth $\chi_u(h)$ estimate for each MCMC iteration. Then, combining across MCMC iterations, we compute the pointwise average curves and their credible bands. The results in Figure \ref{chi_enve} demonstrate that our model provides a sensible approximation to the extremal dependence structure of the mis-specified models. Borrowing strength across locations, the parametric estimators of $\chi_u(h)$ are much more reliable than the nonparametric ones in that they are able to discriminate between the two asymptotic classes via estimating a relatively small number of parameters. Especially for the \citet{huser2017bridging} models, the dependence strength diminishes gradually as $\beta$ becomes greater, eventually resembling the Gaussian copula. With the extra tail flexibility, our model is able to accurately capture these features.

\section{Data Analysis}\label{data_analysis}
We consider daily observations of Fosberg Fire Weather Index (FFWI) from 1974 to 2015 at 93 monitoring stations over parts of the Great Plains, mainly from Central Great Plains to South Texas Plains \citep{dunn2012hadisd}.  Figure \ref{GreatPlains} shows the observation locations as black triangles. The FFWI aims to quantify potential wildfire threat. It is a single number summary calculated from temperature, wind speed, and relative humidity; larger index values signify higher flame lengths and more rapid drying \citep{fosberg1978weather}. Due to human activity and changes in the grassland ecosystem, the Great Plains region is becoming an important high-risk region of large wildfires. According to a in-depth study conducted by \citet{donovan2017surging}, the average total area burned by fires \replace{in Great Plains region between 2005 and 2014 was in the millions of hectares per year. In 2017, 809,380 hectares were lost to wildfires in a single week in Texas, Oklahoma, and Kansas alone \citep{herskovitz2017wildfires}.} As a consequence of the prevailing hot and dry air, wildfires in this region are particularly more concentrated in the spring, feasting on grasses made dry by long-term drought. Modeling the tail behavior of FFWI and studying the extremes of the process could have major implications for wildfire planning.

\begin{figure}
    \centering
    \includegraphics[width=0.75\linewidth]{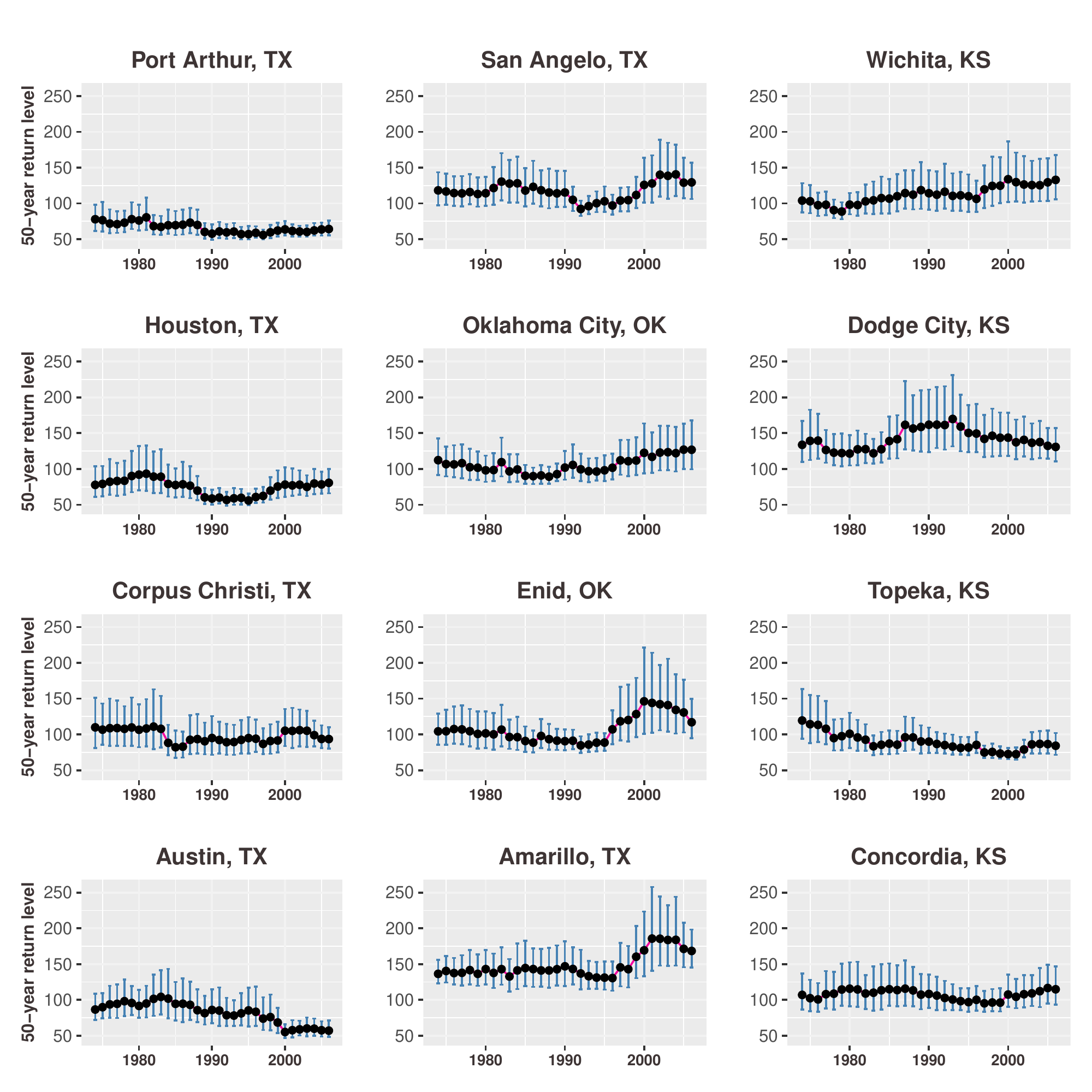}
    \caption{Point estimates and 95\% confidence intervals for 50-year return levels of FFWI at 12 randomly selected stations. For each station and for each year, we estimate the return level using annual maxima in a 10-year sliding window.  There is no evident systematic trend in the return level, so the assumption of marginal parameters that are constant in time is deemed appropriate.}
    \label{Marginal}
\end{figure}

To ensure the independence over time and avoid seasonal effects, we take the maxima of the FFWI values over ten-day intervals during the spring season. Figure \ref{Marginal} shows the 50-year return levels estimated using the block maxima with 10-year sliding windows for 12 randomly-selected stations. There is no clear evidence for a systematic increase or decrease in the return levels. Although treating meteorological data as constant over time is often problematic, particularly for temperature data, the behavior in Figure \ref{Marginal} suggests that an assumption of constant marginal parameters in time is appropriate. 

To account for the physical features of the terrain in the Great Plains, we describe the scale parameter in $\boldsymbol{\theta}_{GPD}$ by the trend surface:
\begin{equation}
    \sigma(\boldsymbol{s})=\beta_0 + \beta_1 \text{lon}(\boldsymbol{s}) + \beta_2 \text{lat}(\boldsymbol{s})
    \label{trendsurf}
\end{equation}
where $\text{lon}(\boldsymbol{s})$ and $\text{lat}(\boldsymbol{s})$ are the longitude and latitude of the stations at which the data are observed.  We constrain the joint prior of $(\beta_0, \beta_1, \beta_2)'$ such that the support of $\sigma(\boldsymbol{s})$ is the positive real line.  We model the shape parameter $\xi(\boldsymbol{s})$ as constant over the spatial domain, as suggested by exploratory analysis (see Appendix \ref{sec:trendsurf}).

Similar to the procedure in Section \ref{mis-specified}, to obtain starting values for MCMC, we fit  generalized Pareto distributions to model events above the 98\% quantile, $u_{98}$, and empirical distributions to those below $u_{98}$, of the time series at each station separately, and then use the fitted models to transform the data to have noisy scale mixture distributions. We then ran the MCMC chain for 50,000 iterations thinned by 10 steps and discarded a burn-in period of 25,000 iterations.

\subsection{Model Evaluation}
First and foremost, we examine whether the estimate $\delta$ falls within $(0,1/2]$ or $(1/2,1)$ in accordance to whether the data-generating process is asymptotically independent or dependent. Table \ref{mcmc_tab} reports the posterior means and 95\% credible intervals for the model parameters. Trace plots for $\delta$ and $\tau$ can be found in Appendix \ref{mcmc_results}. For this dataset, the MCMC results show that the range of the mixing parameter $\delta$ is close to the interface between the two dependence class, while demonstrating asymptotic dependence. Nonetheless, the value of $\delta$ being close to 1/2 means that $\chi_u(h)$ will still decrease with $u$ before eventually reaching a positive limit.
\begin{table}
\begin{adjustbox}{max width=\textwidth}
\begin{tabular}{c|cccc}\hline
 & $\rho$ & $\nu$ & $\tau$ & $\delta$  \\\hline
\textbf{Posterior Mean} & 0.504 & 0.344 & 1.837 & 0.530 \\
\textbf{95\% Credible Interval} & (0.483, 0.524) & (0.324, 0.362) & (0.331, 3.963) & (0.528, 0.531)\\\hline
 & $\beta_0$ & $\beta_1$ & $\beta_2$ & $\xi$ \\\hline
\textbf{Posterior Mean} & -17.531 & -0.284 & -0.1873 & 0.3392 \\
\textbf{95\% Credible Interval} & (-17.614, -17.401) & (-0.285, -0.283) & (-0.1874, -0.1871) & (0.3385, 0.3396)\\
\hline
\end{tabular}
\end{adjustbox}
\caption{Posterior mean and 95\% highest posterior density credible interval.}
\label{mcmc_tab}
\end{table}

\begin{figure}
    \centering
    \includegraphics[width=0.32\linewidth]{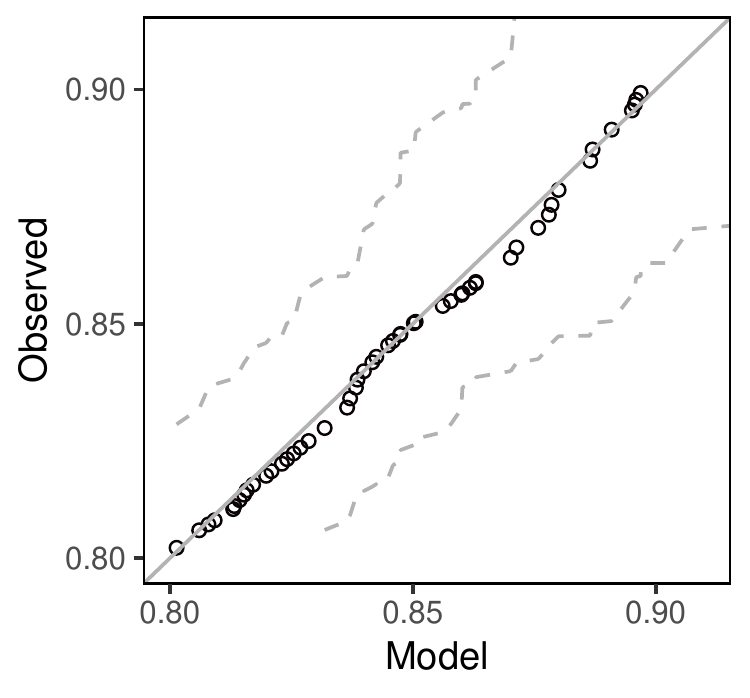}
    \includegraphics[width=0.32\linewidth]{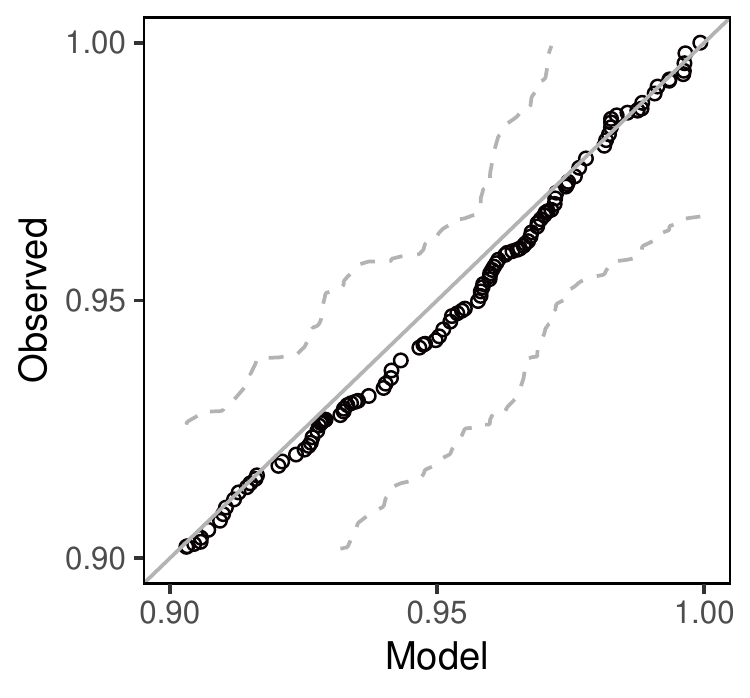}
    \caption{Comparisons of the observed and predicted mean (left) and maxima (right) for 5 locations held out for model validations. Overall 95\% confidence envelopes are also shown. Note that the values are transformed marginally to uniform using the empirical distribution functions at each location.}
    \label{qqplot}
\end{figure}

To better evaluate the model fit, we randomly hold out 5 stations for validation purposes, and exclude them from the MCMC analysis; see the red points in Figure \ref{GreatPlains}. We then compare the empirical distributions of mean, and maxima for each time point at the 5 held-out stations with those simulated with parameters from MCMC samples. Though we modeled our data as censored observations, the model may still work a bit further into the center of the distribution. Results are displayed in Figure \ref{qqplot} where we only show values exceeding 80\% threshold for mean, and 90\% threshold for maxima. We can see that the fit displays a good match against the observed mean and maxima in the upper quantiles. 

\begin{figure}
    \centering
    \includegraphics[width=0.42\linewidth]{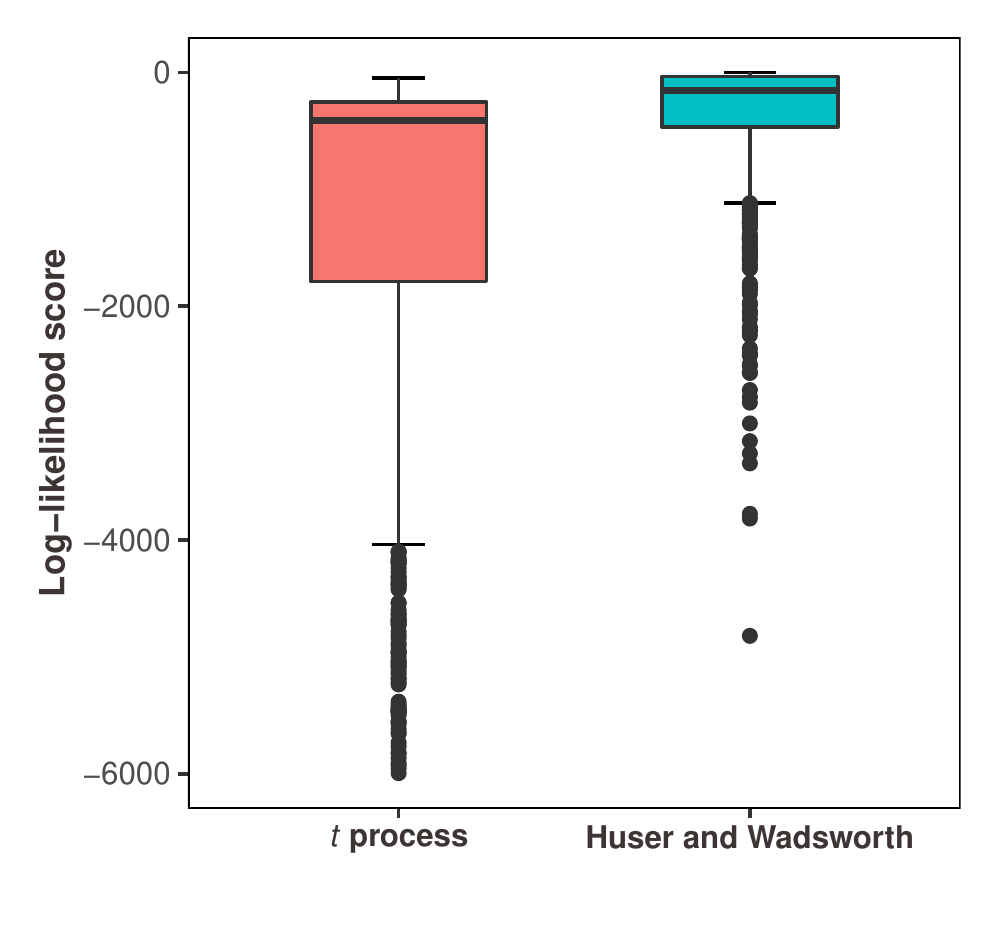}
    \includegraphics[width=0.42\linewidth]{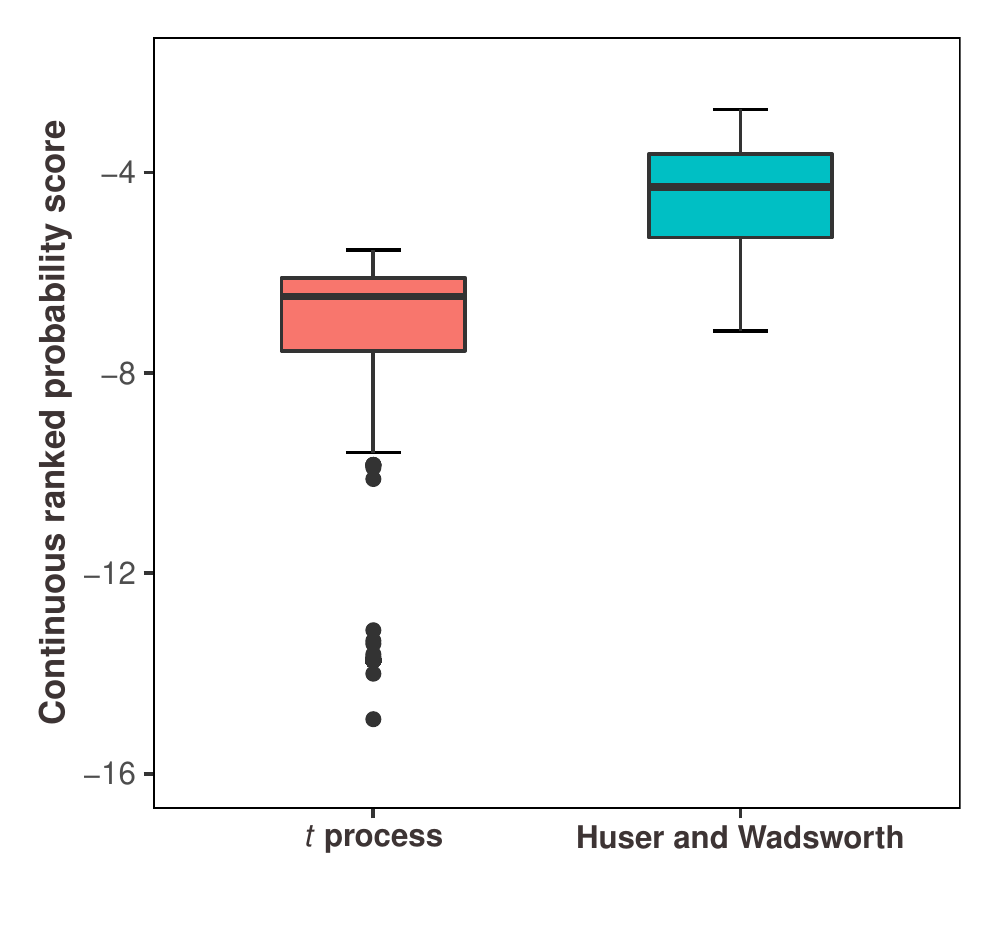}
    \caption{Comparisons of the log-likelihood scores and CRPS between our model and a similar model with $t$ latent process.  In both panels, higher values indicate better model fit.}
    \label{logscore}
\end{figure}

As a comparison, we change $\boldsymbol{X}^*$ to be a $t$ process, and re-fit the model using MCMC. We apply proper scoring rules \citep{gneiting2007strictly}, log scores and continuously-ranked probability scores (CRPS), to compare the quality of probabilistic forecasts. While running MCMC, we interpolate the latent process $\boldsymbol{X}^*$ at the held-out locations for each iteration using the full conditional likelihood. Plugging the predictive draws at the held-out observations into the equation \eqref{fullY}, we obtain the log score (simply the log-likelihood) as
\begin{equation}
    V=\log \left(\prod_{t=1}^T\prod_{j=1}^5 \phi(y_t(\boldsymbol{r}_j) \given \boldsymbol{X}^*_t,\cdots)\right)
    \label{logscores}
\end{equation}
where $\{\boldsymbol{r}_j\given j=1,\ldots,5\}$ are the validation stations. The left panel of Figure \ref{logscore}  compares the log scores between two models, showing that the transformed Gaussian scale mixture process clearly outperforms the $t$ process. Additionally, we calculate the CPRS \citep{matheson1976scoring} for both models,
\begin{equation}
    \text{CRPS}(F_Y,y) = -\int_{-\infty}^{\infty} (F_Y(z)-\mathbbm{1}\{y\leq z\})^2 dz
    \label{CRPS}
\end{equation}
where $F_Y$ is the marginal distribution estimated using parameters using one MCMC iteration, and $y$ is the observed value. The right panel of Figure \ref{logscore} shows the averaged CRPS for the two models, where our model clearly has better results. 

\begin{figure}
    \centering
    \includegraphics[width=0.45\linewidth]{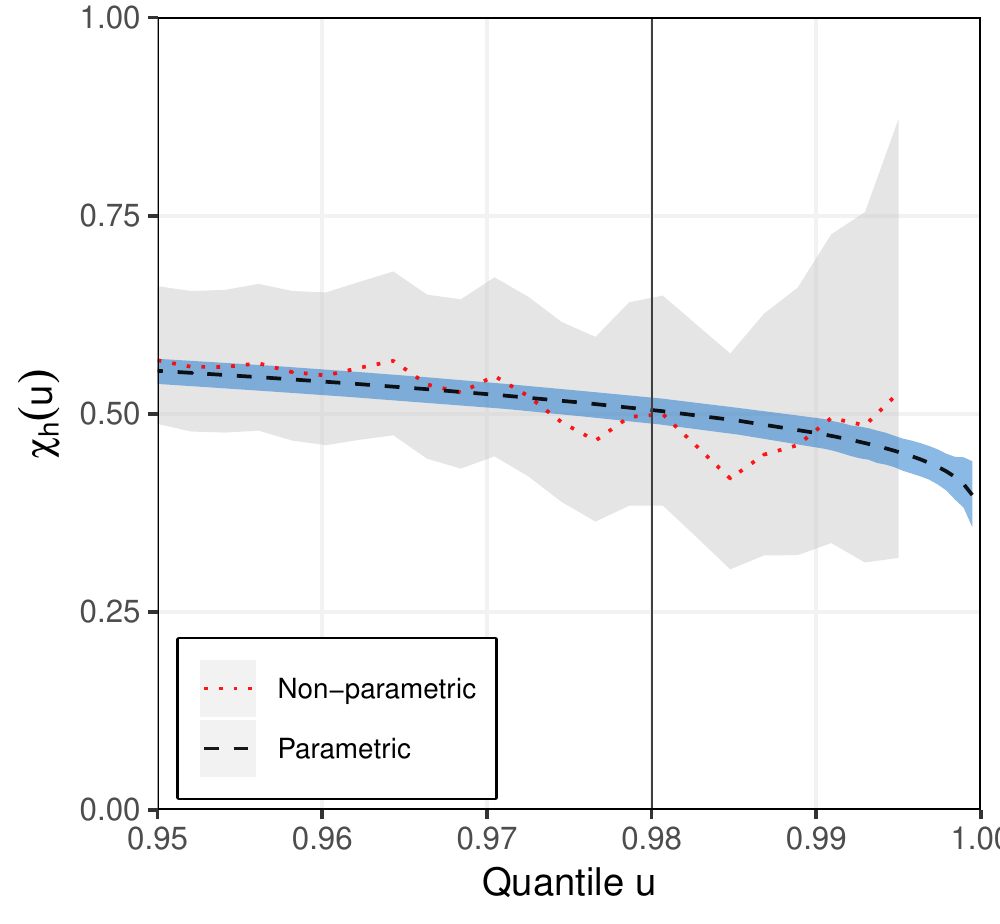}
    \includegraphics[width=0.45\linewidth]{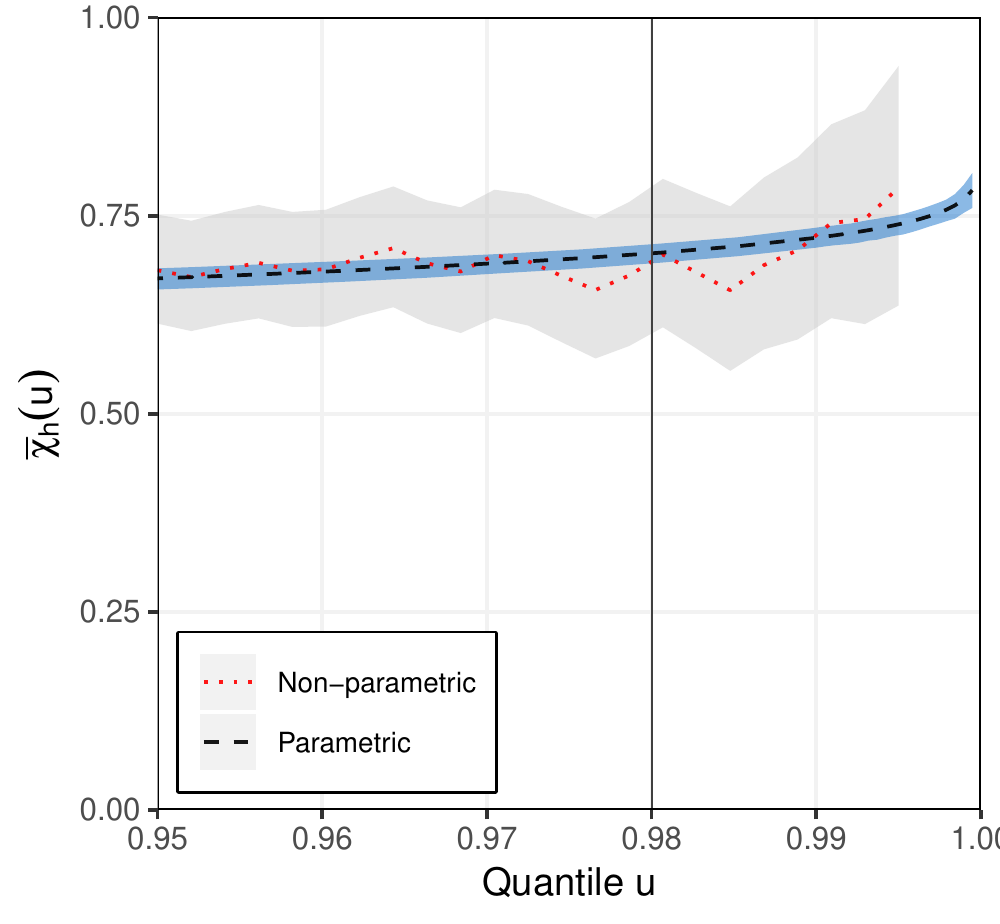}
    \caption{Empirical estimates (dashed red lines and gray envelopes) of $\chi_u(h)$ (left) and $\bar{\chi}_u(h)$ (right) for the FFWI, $u\in [0.95,0.995]$, for two points at distance $\|\boldsymbol{s}_1-\boldsymbol{s}_2\|=20 km$. The blue shaded areas are 95\% credible envelopes obtained from the posterior samples of fitted model, and the black dashed lines show the averaged curve. The vertical line is the threshold used when fitting the dependence model.}
    \label{Chibar}
\end{figure}

Similar to Figure \ref{chi_enve}, we show the empirical and model-based values of $\chi_u(h)$ and $\bar{\chi}_u(h)$ for the block maxima of the FFWI in Figure \ref{Chibar}, which confirms that our model captures the extremal spatial dependence in the data quite well. The quantity $\bar{\chi}_u(h)$  is an alternative dependence measure useful in the situation $\chi(h)=0$, and is defined as
\begin{equation*}
    \bar{\chi}_h(u)=\frac{2\log(1-u)}{\log P(F_j(X_j)>u, F_k(X_k)>u)}-1 \to 2\eta\editLZ{_X}(h)-1, \qquad u \to 1.
\end{equation*}
Recall that the posterior mean for $\delta$ is greater than 0.5, which means $\bar{\chi}_h(u)\to 1$ as $u\to 1$. Interestingly, the black dashed curve of the right panel of Figure \ref{Chibar} seems to have a limit less than 1. This is because to attain the correct limit in this case, we would need to compute $\bar{\chi}_h(u)$ for values of $u$ that are  \textit{very} close to 1, which is very difficult numerically.

\subsection{Results}
 
\begin{figure}
    \centering
    \includegraphics[width=0.45\linewidth]{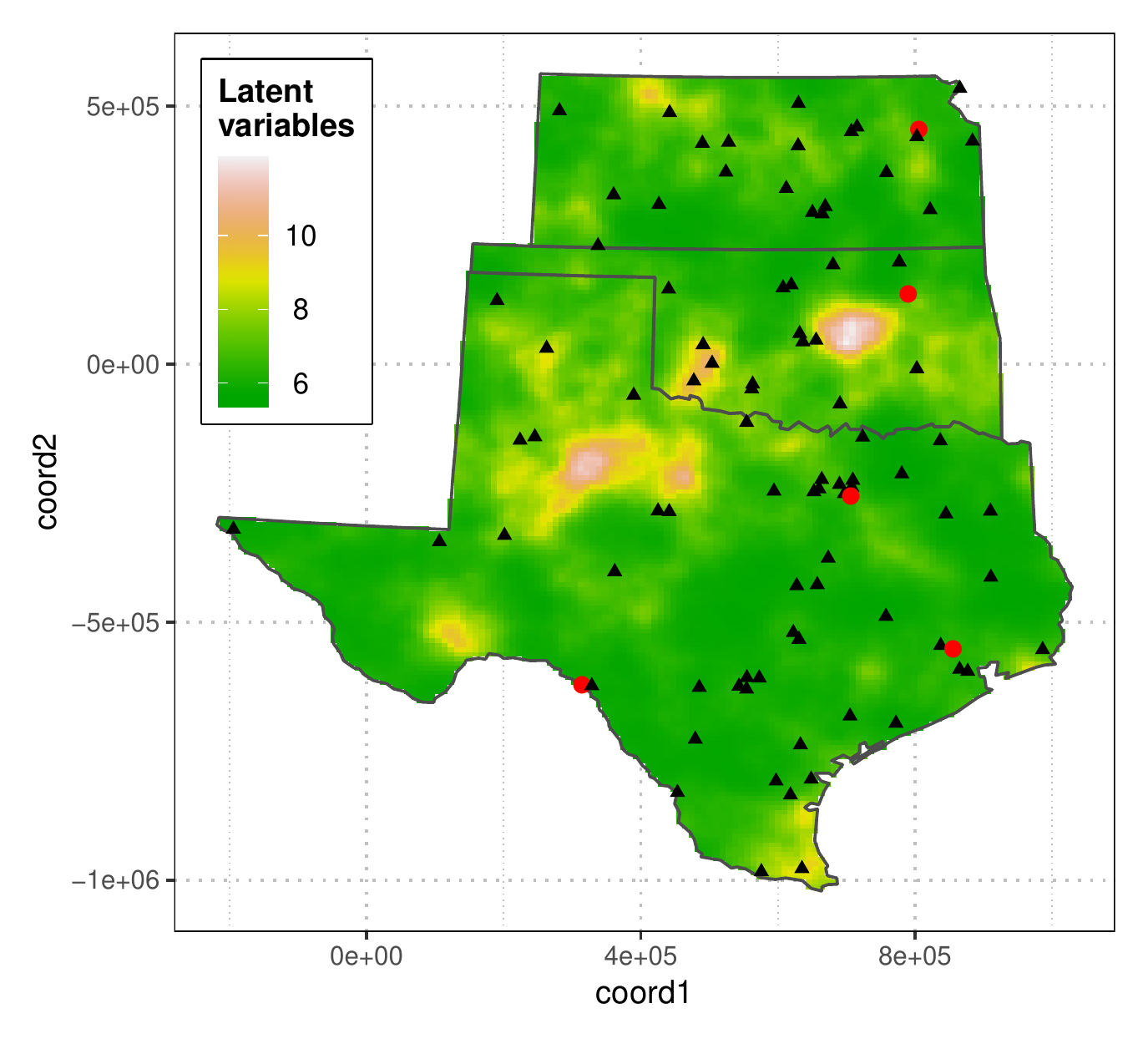}
    \includegraphics[width=0.45\linewidth]{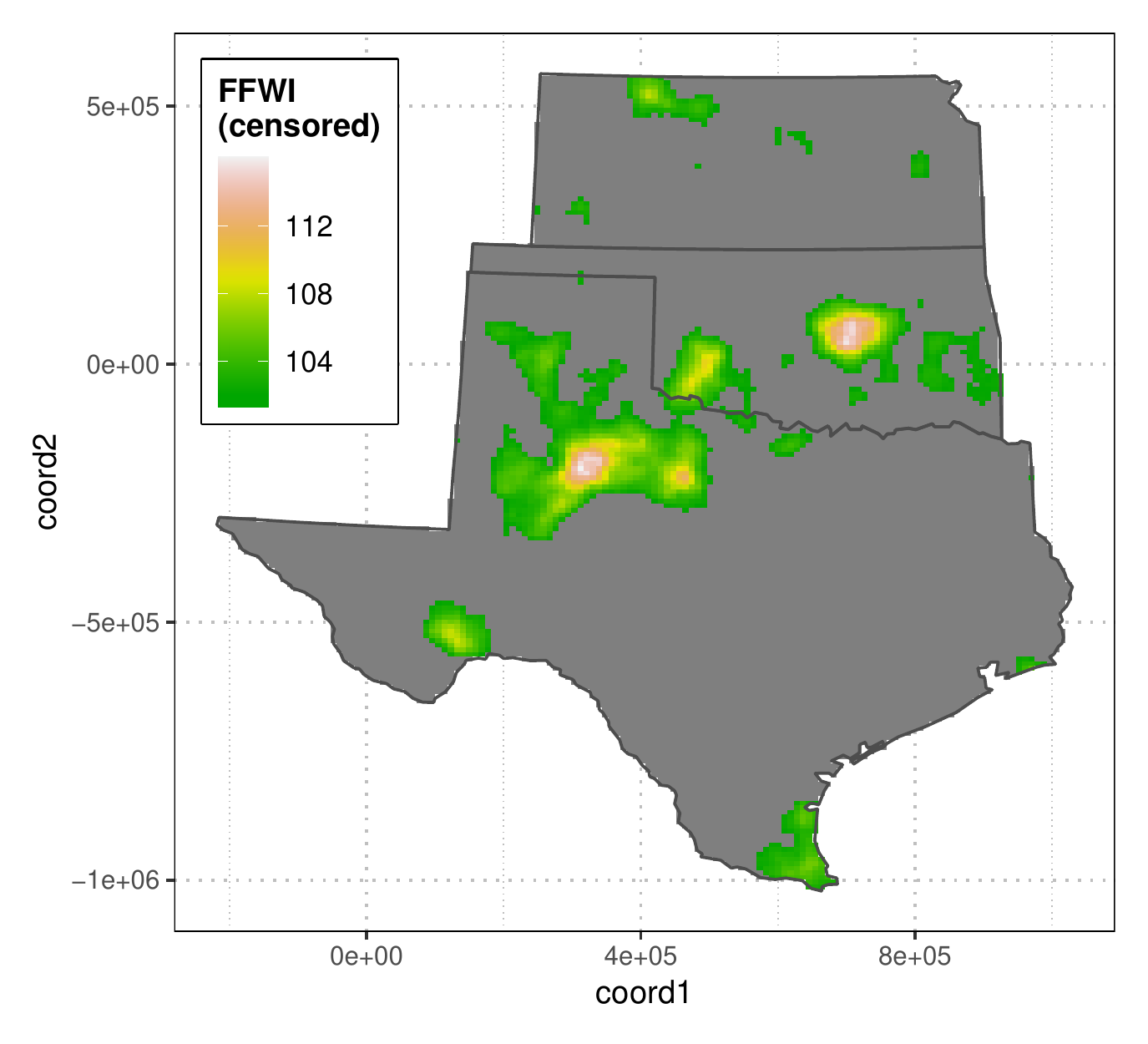}
    \caption{A realization of the process generated from parameters from one MCMC iteration. Left panel shows the latent scale mixture process $X(s)$ \replace{in log scale to better visualize spatial patterns}, with points showing the stations of the 93 gauges. For model checking, the 88 stations marked by triangles were used to fit the models, and those marked with circles were used to validate the models. Right panel shows the same realization, marginally transformed to the scale of the original FFWI data, and then censored. A projected coordinate reference system, NAD83(HARN), is used here so that the axes are in units of meters.}
    \label{GreatPlains}
\end{figure}

To get an idea of what a realization of the fitted process looks like, the left panel of Figure \ref{GreatPlains} shows one realization of the latent $X(s)$ scale mixture process using parameters from one MCMC iteration. The extreme values are mainly concentrated in two small regions. The right panel shows the same realization, now marginally transformed to the scale of the FFWI values. Since we modeled our data as partially censored observations, the map here only displays the areas where threshold exceedances are observed. 

\begin{figure}
    \centering
    \includegraphics[width=0.45\linewidth]{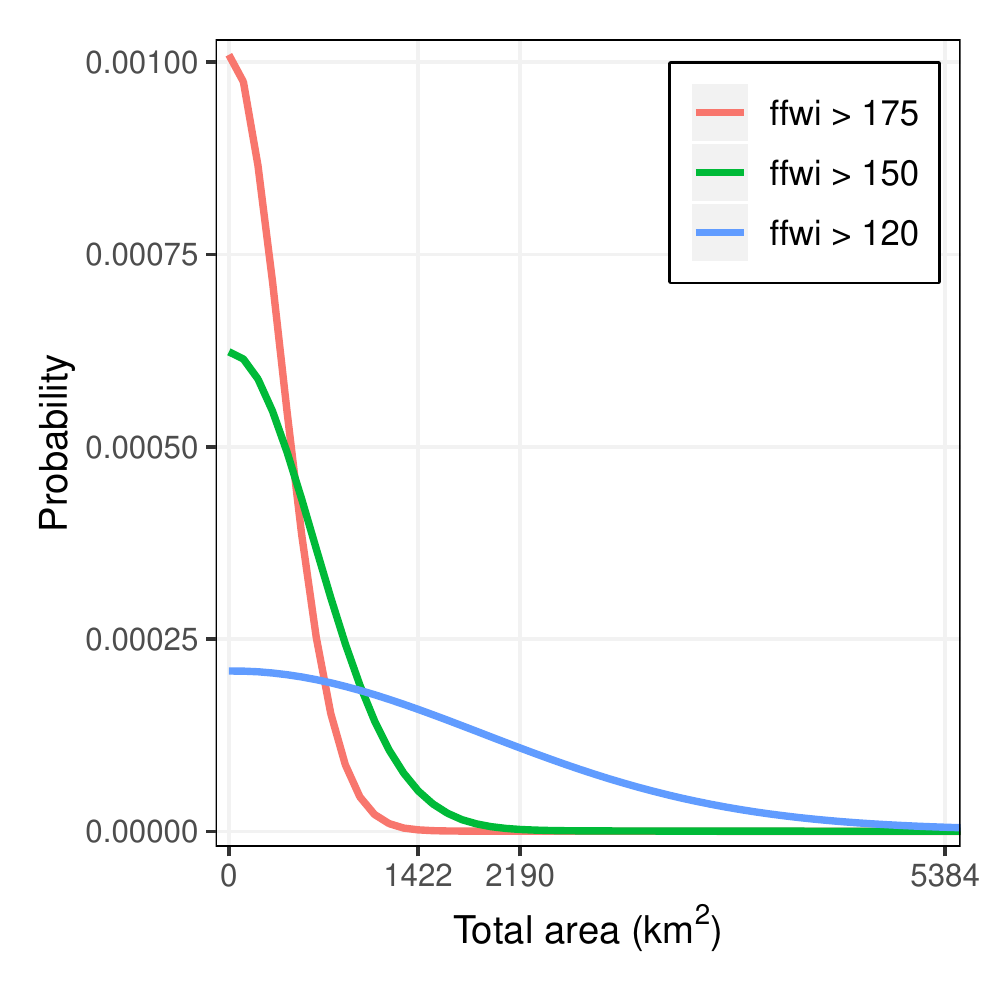}
    \includegraphics[width=0.45\linewidth]{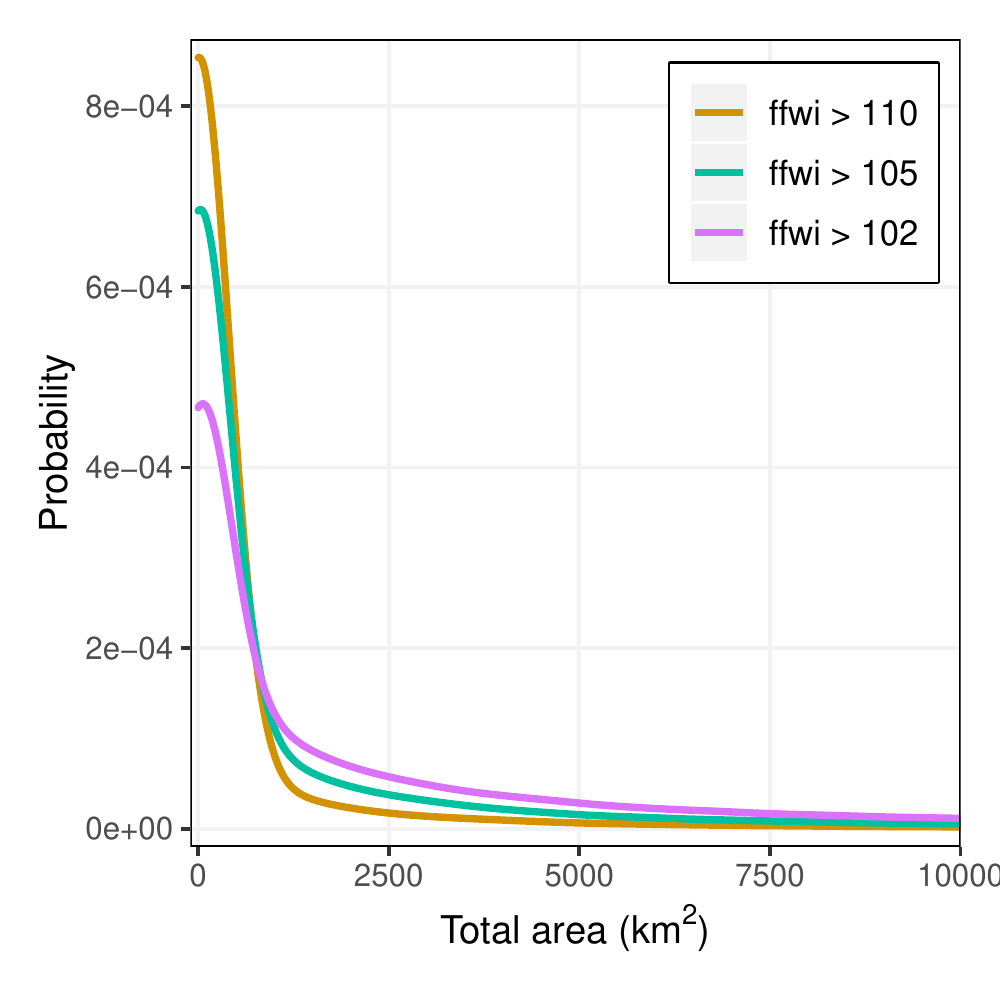}
    \caption{The distribution of the total area that is over a certain threshold. For each MCMC iteration, we simulate 100 processes on a $15km\times 15km$ grid, and count the number of threshold exceedances. The density of the total area is then estimated using Gaussian kernels.}
    \label{Exceedance}
\end{figure}

A quantity of great interest is areal exceedance probabilities, which represent the amount of territory simultaneously at extreme risk for wildfire.  To obtain a Monte Carlo estimate of these joint probabilities, we use parameters from each MCMC sample to simulate 100 processes (on a $15km \times15km$ grid), and calculate the total area that has FFWI over a designated threshold. Figure \ref{Exceedance} shows the results.  These curves represent total area at risk for various FFWI thresholds.  The curves for the higher thresholds decay faster than those for the lower thresholds, which confirms that extreme events simultaneously occurring across large areas becomes less common when the threshold increases. This also shows that the joint tail of the fire threat index exhibits a weakening dependence structure, with more extreme events being more localized.  It is not possible to capture this behavior using limiting extreme value models like max-stable or generalized Pareto processes.

\section{Discussion}\label{sec:discussion}

In this paper, we have proposed a new modeling approach, based on the class of transformed Gaussian scale mixture models, which includes those recently proposed by \citet{huser2017bridging} and \citet{huser2017modeling}. We added a measurement error to the mixture, hence avoiding the need to calculate the onerous $|\mathcal{C}|$-dimensional Gaussian distribution function when dealing with the censored likelihood. We also circumvent the need to draw from a high-dimensional truncated distribution by treating the smooth process as latent and updating repeatedly using MCMC. In addition to its computational advantages, the presence of a measurement error term may also make the model more realistic for data collected by real-world instruments. Indeed, nugget effects are ubiquitous in environmental statistics, not just to represent measurements errors, but also as a result of small scale effects that are not included in the large scale model for spatial dependence.

Even with the presence of the measurement error, the model is still able to capture qualitatively different types of sub-asymptotic dependence behavior of spatial processes. In the case of our modification of the \citet{huser2017modeling} model, a smooth transition between both extremal dependence paradigms takes place in the interior of the parameter space, which enables inference of the dependence class in a simple manner. We proved that all the appealing asymptotic properties found in the original smooth process are preserved in the modified model, regardless of the size of the measurement error variance.

The model allows inference on spatial extreme-value datasets with relatively large numbers of locations. The computational limitations are similar to those of conventional spatial Gaussian process models.  We have defined the model conditionally as a Bayesian hierarchical model, for which standard MCMC techniques can be used to fit the data. Computation is facilitated greatly by parallelizing over time $t$ and migrating some basic linear algebra to \texttt{C/C++} via \texttt{Rcpp}.  Even so, the lack of closed form marginal transformations creates a significant (though embarrassingly parallel) computational challenge that scales with the total number of exceedances, rather than the usual case of scaling with the number of spatial locations.

Despite easing computational limitations associated with the \citet{huser2017bridging} and \citet{huser2017modeling} models, our modified versions inherit the same theoretical limitations. Neither model is able to account for the possibility of independence between observations as the distance between sites becomes large, nor are they able to transition from asymptotic dependence at short range to asymptotic independence at longer range. \citet{wadsworth2019higher} presents an alternative approach to modeling spatial extremes that begins to address these issues. 

\replace{Another interesting possibility to explore would be to include the nugget term inside the link function $g$. This would result in closed-form marginal distributions in some cases, perhaps making computations easier.  However, it would change the dependence structure in ways that would not vanish in the limit, which is not the behavior we were aiming for here, but could be useful nonetheless.}

\if 0\blind{
\section*{Acknowledgements}
We gratefully acknowledge support from NSF grant DMS-1752280 and EPSRC grant EP/P002838/1, along with seed grants from the Institute for CyberScience and the Institute for Energy and the Environment at Pennsylvania State University.  Computations for this research were performed on the Pennsylvania State University’s Institute for CyberScience Advanced CyberInfrastructure (ICS-ACI).
}\fi

\appendix
\section{Technical appendix}
\label{app:proofs}
\subsection{Proof of Proposition~\ref{prop:dep}}
For the proof of Proposition~\ref{prop:dep}, we begin by recalling useful results from the literature. 

The first is Breiman's lemma, see e.g.\ \citet{brieman-1965a} and \citet{cline-1994a}, and a corollary for sums of a regularly varying and light-tailed random variables. 
\begin{lemma}[Breiman]
 \label{lem:breiman}
 Suppose that $Q=ST$ where $\Pr(S>s) \in \RV_{-\alpha}$, $\alpha\geq 0$ and $\E(T^{\alpha+\delta})<\infty$ for some $\delta>0$. Then
 \[
  \Pr(Q>x) \sim \E(T^\alpha)\Pr(S>x),\qquad x\to\infty.
 \]
\end{lemma}

\begin{corollary}
 \label{cor:breiman}
 Suppose that $\P(X>x) \in \RV_{-\alpha}$, $\alpha\in(0,\infty)$, and $\epsilon$ is a random variable such that $\E(e^{\delta \epsilon})<\infty$ for some $\delta>0$. Then
 \[
  \Pr(X+\epsilon>x) \sim \Pr(X>x), \qquad x\to\infty.
 \]
\end{corollary}
\begin{proof}
 Since $\P(X>x) =:\bar{F}_X(x) \in \RV_{-\alpha}$, for positive finite $\alpha$, $\P(e^X>x) = \bar{F}_X(\log(x)) \in \RV_0$. By assumption $\E(e^{\delta \epsilon})<\infty$, and so applying Breiman's Lemma to $e^X e^\epsilon$ yields 
 \[
  \Pr(e^X e^\epsilon > x) \sim \Pr(e^X>x),\qquad x\to\infty,
 \]
from which the result follows.
\end{proof}

The second result relates to sums of Weibull-tailed variables, i.e., those with survival functions satisfying~\eqref{eq:lightsf} and the associated density 
\begin{align}
 f_{X^*}(x) \sim v(x) \exp(-\theta x^\alpha), \qquad v(x) = u(x)(\theta\alpha x^{\alpha-1}) \in \RV_{\kappa+\alpha-1}. \label{eqn:lightdens}
\end{align}
The following lemma can be verified directly from Theorem~4.1 of \citet{asmussen-2017a}, by identifying that the conditions in Section 4 of that paper hold for this subclass when $\alpha>1$.

\begin{lemma}[\cite{asmussen-2017a,balkema-1993}]
\label{lem:weibullsum}
 Let $Y_1, Y_2$ be variables with density~\eqref{eqn:lightdens}, with regularly varying functions $v_j, u_j$ and parameters $\theta_j>0, \alpha_j>1$, $j=1,2$. Then the density and survival function of the convolution satisfy
 \begin{align*}
  f_{Y_1+Y_2}(x) &\sim v_+(x)\exp\{-\psi_{+}(x)\}, & \P(Y_1+Y_2>x) &\sim u_{+}(x)\exp\{-\psi_{+}(x)\},
 \end{align*}
where $\psi_+(x) = \theta_1 q_1(x)^{\alpha_1}+ \theta_2 q_2(x)^{\alpha_2}$, with $q_1(x),q_2(x)$ determined by solving
\begin{align}
 q_1+q_2 = x, \qquad \theta_1\alpha_1 q_1^{\alpha_1-1} = \theta_2\alpha_2 q_2^{\alpha_2-1}, \label{eq:q1q2}
\end{align}
and
\begin{align*}
 v_{+}(x) = \left(\frac{2\pi \psi_{+}''(x)}{\theta_1\alpha_1 (\alpha_1-1)q_1^{\alpha_1-2} \theta_2\alpha_2 (\alpha_2-1)q_2^{\alpha_2-2}}\right)^{1/2} v_1(q_1) v_2(q_2); \qquad u_+(x) =v_+(x)/\psi'_+(x).
\end{align*} 
\end{lemma}
Finally we note also the following two useful inequalities
 \begin{align}
  \Pr(X^*_1+\min(\epsilon_1,\epsilon_2) >x, X^*_2+\min(\epsilon_1,\epsilon_2) >x) \leq &  \Pr(X^*_1+\epsilon_1 >x, X^*_2+ \epsilon_2 >x)  \leq \notag \\
  &\Pr(X^*_1+\max(\epsilon_1,\epsilon_2) >x, X^*_2+\max(\epsilon_1,\epsilon_2) >x), \label{eq:ineq}
 \end{align}
  \begin{align}
  \P(\epsilon_1+\min(X^*_1,X^*_2) >x, \epsilon_2+\min(X^*_1,X^*_2) >x) \leq   &\P(X^*_1+\epsilon_1 >x, X^*_2+ \epsilon_2 >x)  \leq \notag \\
  &\P(\epsilon_1+\max(X^*_1,X^*_2) >x, \epsilon_2+\max(X^*_1,X^*_2) >x) \label{eq:ineq2},
 \end{align}
where the two lower bounds are equal. We can now prove Proposition~\ref{prop:dep}. 

\begin{proof}[Proof of Proposition~\ref{prop:dep}]
1) When $X^*$ has a regularly varying tail, Corollary~\ref{cor:breiman} provides that $\P(X^*+\epsilon>x) \sim \P(X^*>x)$. For the existence of $\eta_{X^*}$, the variable $\min(X_1^*, X_2^*)$ also has a regularly varying tail. Further,
\[
\E(e^{\delta\max(\epsilon_1,\epsilon_2)}) = \E\left[\max(e^{\delta \epsilon_1},e^{\delta \epsilon_2})\right] \leq \E\left[e^{\delta \epsilon_1}+e^{\delta \epsilon_2}\right] = 2e^{\delta^2\sigma^2/2}<\infty,
\]
so Corollary~\ref{cor:breiman} thus gives $\Pr(X^*_1+\max(\epsilon_1,\epsilon_2) >x, X^*_2+\max(\epsilon_1,\epsilon_2) >x) \sim \Pr(X^*_1>x,X^*_2>x)$, and similarly for the lower bound. Hence $\Pr(X^*_1+\epsilon_1 >x, X^*_2+ \epsilon_2 >x) \sim \Pr(X^*_1>x, X^*_2>x)$. Therefore it follows that $\chi_{X^*}=\chi_X$. If $f(x)\sim g(x)$ for $f(x) \to 0$ then $\log f(x) \sim \log g(x)$, so the result for $\eta$ follows as well. 

When $X^*$ has a Weibull-like tail with Weibull index $\alpha<1$, then $\P(e^{X^*}>x) \sim u(\log x)\exp\{-\theta(\log x)^\alpha\} \in \RV_0$, and the rest of the argument follows as above.

\noindent
2) When $\alpha=1$, $\P(e^{X^*}>x) \sim u(\log x) x^{-\theta}\in \RV_{-\theta}$. Breiman's lemma now yields
\begin{align}
 \P(e^{X^*}e^\epsilon>x) \sim \E(e^{\theta \epsilon}) \P(e^{X^*}>x). \label{eqn:a1marg}
\end{align}
As a consequence, inequality~\eqref{eq:ineq} does not lead to a precise asymptotic relationship for $\P(X^*_1+\epsilon_1 >x, X^*_2+ \epsilon_2 >x)$, but rather that it is asymptotically bounded within the range
\begin{align}
  \left[\E\left[e^{\theta\min(\epsilon_1,\epsilon_2)}\right],  \E\left[e^{\theta\max(\epsilon_1,\epsilon_2)}\right]\right] \P(\min(X_1^*,X_2^*) >x). \label{eqn:a1jt}
\end{align}
Combining~\eqref{eqn:a1marg} and~\eqref{eqn:a1jt}, we get the stated bound for $\chi_{X^*}$. It follows also that $\log \P(X>x) \sim \log\P(X^*>x)$ and $\log\P(X_1>x,X_2>x) \sim \log\P(X^*_1>x,X^*_2>x)$, and hence $\eta_{X^*} = \eta_{X}$.

\noindent
3) Here we use Lemma~\ref{lem:weibullsum}, where different values of $q_1,q_2$ are found for the three cases $\alpha \in(1,2)$, $\alpha=2$ and $\alpha>2$. To make notation more obvious, we replace $q_1, q_2$ with $q_*, q_\epsilon$ for summation of $X^*, \epsilon$, and $q_*^\wedge, q_\epsilon^\vee$ etc., for summation of $\min(X_1^*,X_2^*)$ and $\max(\epsilon_1,\epsilon_2)$, for example. The three cases are considered for the value of $\alpha_*$; we always have $\alpha_{\epsilon}=2$.

\noindent
a) For $\alpha_* \in (1,2)$,
\begin{align*}
  q_*(x) &= x - \frac{ \alpha_*\theta_*}{2 \theta_\epsilon}x^{\alpha_*-1}[1+o(1)],\\
  q_\epsilon(x) &= \frac{ \alpha_*\theta_*}{2 \theta_\epsilon}x^{\alpha_*-1}[1+o(1)].
 \end{align*}
Therefore
\begin{align}
 \psi_+(x) &= \theta_* q_*(x)^{\alpha_*}+ \theta_\epsilon q_\epsilon(x)^{2} \sim \theta_* x^{\alpha_*}, \label{eq:psi+}
\end{align}
which implies $-\log\P(X^*+\epsilon>x) \sim -\log\P(X^*>x)$.

To understand the joint behavior, we again use~\eqref{eq:ineq}. Note that
\begin{align*}
 \P(\max(\epsilon_1,\epsilon_2)>x) =1-\Phi(x/\tau)^2 \sim 2 \phi(x/\tau)/(x/\tau) = u_\epsilon^\vee(x) e^{-x^2/(2\tau^2)}. \\
 \P(\min(\epsilon_1,\epsilon_2)>x) =[1-\Phi(x/\tau)]^2 \sim \phi(x/\tau)^2/(x/\tau)^2 = u_\epsilon^\wedge(x) e^{-x^2/\tau^2}.
\end{align*}
Consequently, $\min(\epsilon_1,\epsilon_2)$ and $\max(\epsilon_1,\epsilon_2)$ both have $\alpha=2$, with different $\theta$. As $\theta_\epsilon$ does not affect the leading order behavior of $\psi_+$ in~\eqref{eq:psi+}, the tails of $\min(X^*_1,X^*_2) + \min(\epsilon_1,\epsilon_2)$ and $\min(X^*_1,X^*_2) + \max(\epsilon_1,\epsilon_2)$ both have $\psi_+(x) \sim \theta_*^\wedge x^{\alpha_*}$ in the exponent, where 
\begin{align*}
 \P(\min(X^*_1,X^*_2) > x) \sim u_{*}^\wedge(x) \exp\{-\theta_*^\wedge x^{\alpha_*}\}.
\end{align*}
Consequently $-\log\P(X^*_1+\epsilon_1 >x, X^*_2+ \epsilon_2 >x) \sim -\log \P(X^*_1 >x, X^*_2 >x)$ and so $\eta_{X^*}=\eta_X$.

\noindent
b) When $\alpha_*=\alpha_\epsilon=2$, solving equations~\eqref{eq:q1q2} provides
\begin{align*}
 q_*(x) = \left(\frac{\theta_\epsilon}{\theta_*+\theta_\epsilon}\right)x, \qquad  q_\epsilon(x) = \left(\frac{\theta_*}{\theta_*+\theta_\epsilon}\right)x,
\end{align*}
and 
\begin{align*}
 \psi_+(x) =  \left(\frac{\theta_\epsilon\theta_*}{\theta_*+\theta_\epsilon}\right) x^2.
\end{align*}
For the margins, and maximum $\theta_\epsilon = \theta_\epsilon^\vee = 1/(2\tau^2)$, whilst for the minimum, $\theta_{\epsilon}^\wedge = 1/\tau^2$. We can now use both sets of inequalities~\eqref{eq:ineq} and~\eqref{eq:ineq2}, to give that the range of $\eta_{X}$ is 
\begin{align}
\left[\frac{\theta_*}{2\theta_*^{\wedge}}\frac{\theta_*^{\wedge}+1/\tau^2}{\theta_*+1/(2\tau^2)}, \min\left\{\frac{\theta_*}{\theta_*^{\wedge}} \frac{\theta_*^{\wedge}+1/(2\tau^2)}{\theta_*+1/(2\tau^2)},\frac{\theta_*}{2\theta_*^{\vee}}\frac{\theta_*^\vee + 1/\tau^2}{\theta_* + 1/(2\tau^2)}\right\}\right]. \label{eqn:a2etarange}
\end{align}
Some simplification arises upon noting that $\theta_*^\vee = \theta_*$, since 
\begin{align*}
\P(\max(X_1^*,X_2^*)>x) &= 2\P(X_1^*>x) -\P(X_1^*>x,X_2^*>x) \\
&=2\P(X_1^*>x)\left\{1-\P(X_2^*>x|X_1^*>x)/2\right\}\\
&\sim 2\P(X_1^*>x)\{1-\chi_{X^*}/2\}, \qquad x \to \infty.
\end{align*}
As $\tau^2\to 0$, i.e., as the nugget effect disappears, the endpoints converge to $\theta_*/\theta_*^\wedge = \eta_{X^*}$. As $\tau^2 \to \infty$, i.e., as the nugget effect dominates, both endpoints converge to $1/2$.

\noindent
c) When $\alpha_*>\alpha_\epsilon=2$, solving equations~\eqref{eq:q1q2} provides
\begin{align*}
 q_*(x) &= \left(\frac{2\theta_\epsilon}{\theta_*\alpha_*} x\right)^{1/(\alpha_*-1)}[1+o(1)],\\
 q_\epsilon(x) &= x- q_*(x) \sim x,
\end{align*}
and $\psi_+(x) \sim \theta_\epsilon x^2$. Using inequality~\eqref{eq:ineq2}, 
\begin{align*}
 -\log \P(X^*_1+\epsilon_1 >x, X^*_2+ \epsilon_2 >x) \sim \theta_\epsilon^{\wedge} x^2;
\end{align*}
noting that $\theta_\epsilon = 1/(2\tau^2)$ and $\theta_{\epsilon}^\wedge = 1/\tau^2$ leads to $\eta_{X}=1/2$.

\end{proof}


\subsection{Hierarchical model}\label{hierarch_mod}
Under spatiotemporal setting, we observe $\{Y_t(\boldsymbol{s})\given t=1,\cdots, T, \boldsymbol{s}\in \mathcal{S}\}$, and the temporal dependence is ignored. Then the hierarchical model can be described as
\begin{equation}\label{hierarchicalModel}
    \begin{split}
        Y_t(\boldsymbol{s})\given \boldsymbol{X}_t^*,\tau^2,\boldsymbol{\theta}_R,\boldsymbol{\theta}_{GPD},p&=T^{-1}(X_t^*(\boldsymbol{s})+\epsilon(\boldsymbol{s})),\\
        \boldsymbol{X}_t^*\given R_t,\boldsymbol{\theta}_C&=R_t\cdot g(\boldsymbol{Z}_t),\\
        R_t\given \boldsymbol{\theta}_R&\sim F_R,\\
        \tau^2&\sim IG(\alpha,\beta),\\
        \boldsymbol{\theta}_{GPD}=(\sigma,\xi)&\sim \text{halfCauchy}(1)\cdot U(\xi;-0.5,0.5),\\
        \boldsymbol{\theta}_C=(\rho,\nu)&\sim \text{halfCauchy}(1)\cdot\text{halfCauchy}(1),
    \end{split}
\end{equation}
where $\boldsymbol{\theta}_R$ corresponds to $\delta\in [0,1]$ for the \citet{huser2017modeling} model and $\beta\in [0,\infty)$ for the \citet{huser2017bridging} model, and $Z_t$ is a Gaussian process with Mat\'{e}rn covariance $\mathbf{\Sigma}_{\theta_C}$. The priors for $\boldsymbol{\theta}_R$ are $\delta\sim U(0,1)$ and $\beta\sim \text{halfCauchy}(1)$ respectively.

In Section \ref{sec:bayes}, we have formulated the full conditional likelihood $\varphi(Y_t(\boldsymbol{s})\given\boldsymbol{X}_t^*,\tau^2,\boldsymbol{\theta}_R,\boldsymbol{\theta}_{GPD},p)$ for a fixed time and location; see Equation \eqref{fullY}. Next we are going to work out the conditional joint likelihood $\varphi(\boldsymbol{X}^*_t\given R_t, \boldsymbol{\theta}_C)=\varphi(X^*_t(\boldsymbol{s}_1),\cdots, X^*_t(\boldsymbol{s}_D)\given R_t, \boldsymbol{\theta}_C)$ for a fixed time $t$. Keeping the notations same as the main text, we denote $g(\cdot)$ as the link function that transforms the scale of the latent Gaussian process. For the model of \citet{huser2017bridging}, $g(\cdot)$ is an identity function, and thus $\boldsymbol{X}_t^*$ conditioning on $R_t$ remains to be a multivariate Gaussian random variable. For the \citet{huser2017modeling} model, the conditional density needs more careful scrutiny.

\begin{lemma}\label{conditionalR}
Let $\boldsymbol{X}^*=(X^*_1,\cdots,X^*_D)=(R\cdot g(Z_1),\cdots,R\cdot g(Z_D))$, where $(Z_1,\cdots, Z_D)\sim N(0,\mathbf{\Sigma}_{\boldsymbol{\theta}_C})$. For the \citet{huser2017modeling} model, the density conditional on $R$ is
\begin{equation*}
    \varphi(x^*_1,\cdots,x^*_D\given R,\boldsymbol{\theta}_C)=\frac{1}{\sqrt{|\mathbf{\Sigma}_{\boldsymbol{\theta}_C}|}}\exp\left\{-\frac{1}{2}g^{-1}\left(\frac{\boldsymbol{x}^*}{R}\right)^T (\Sigma^{-1}_{\theta_C}-I_D)g^{-1}\left(\frac{\boldsymbol{x}^*}{R}\right)\right\} \cdot\frac{R^D}{(\prod_{i=1}^D x^*_i)^2},
\end{equation*}
where $g(\cdot)=1/\{1-\Phi(\cdot)\}$.
\end{lemma}
{\it Proof. }
First note that 
\[
(Z_1,\cdots,Z_D)=(g^{-1}\left(\frac{X^*_1}{R}\right),\cdots,g^{-1}\left(\frac{X^*_D}{R}\right))=(\Phi^{-1}\left(1-\frac{R}{X^*_1}\right),\cdots,\Phi^{-1}\left(1-\frac{R}{X^*_D}\right)).
\]
Fixing $R$ as a constant, we can apply the chain rule to obtain
\begin{equation}\label{eq3}
\begin{split}
 \evalat[\Big]{\frac{d Z_1}{d X^*_1}}{X^*_1=x} &=\sqrt{2\pi}\exp\left\{\frac{[\Phi^{-1}(1-R/x)]^2}{2}\right\}\cdot \frac{R}{x^2}\\
 &=\sqrt{2\pi}\exp\left\{\frac{g^{-2}(x/R)}{2}\right\}\cdot \frac{R}{x^2}.
\end{split}
\end{equation}
Then the Jacobian matrix can be written as
\[
\frac{d \boldsymbol{Z}}{d \boldsymbol{X^*}} =
\begin{bmatrix}
\frac{d Z_1}{d X^*_1} &  & \\
 & \ddots &  \\
 &  & \frac{d Z_D}{d X^*_D}
\end{bmatrix},
\]
and 
\[\det\left(\evalat[\Big]{\frac{d \boldsymbol{Z}}{d \boldsymbol{X^*}}}{\boldsymbol{X}^*=\boldsymbol{x}^*}\right)\stackrel{\eqref{eq3}}{=}(\sqrt{2\pi})^D\exp\left\{\frac{1}{2}\sum_{i=1}^Dg^{-2}(x^*_i/R)\right\} \cdot\frac{R^D}{(\prod_{i=1}^D x^*_i)^2}\]

Denote $f_Z$ as the density function of $N(0,\mathbf{\Sigma}_{\boldsymbol{\theta}_C})$. Then 
\[
f_Z(g^{-1}(x_1^*/R),\cdots, g^{-1}(x_D^*/R))=\frac{1}{(2\pi)^{D/2}\sqrt{|\mathbf{\Sigma}_{\boldsymbol{\theta}_C}|}}\exp\left\{-\frac{1}{2}g^{-1}\left(\frac{\boldsymbol{x}^*}{R}\right)^T \mathbf{\Sigma}^{-1}_{\boldsymbol{\theta}_C}\;g^{-1}\left(\frac{\boldsymbol{x}^*}{R}\right)\right\}.
\]
Apply the change of variables formula
$$\varphi(\boldsymbol{x}^*\given R, \boldsymbol{\theta}_C)=f_Z(g^{-1}(\boldsymbol{x}^*/R))\left|det\left(\frac{d \boldsymbol{Z}}{d \boldsymbol{X^*}}\right)\right|,$$
and we will get the desired result.
\hfill\(\Box\)

\vskip 0.3cm
When updating certain variables for a random walk Metropolis-Hastings step, we calculate the corresponding full conditional likelihood for both current values and proposed values. Since $\varphi(Y_t(\boldsymbol{s})\given\boldsymbol{X}_t^*,\tau^2,\boldsymbol{\theta}_R,\boldsymbol{\theta}_{GPD},p)$ entails most parameters and latent variables, we will evaluate this density over and over again. Looking back at its analytic form in \eqref{fullY}, the most compute-intensive part is to calculate the marginal transformation $T(\cdot)$ as defined in \eqref{marginal_Transform}, which requires the computation of the marginal quantile function of the noisy process, i.e. $F_{X\given\btheta_R,\tau^2}^{-1}$. Note that the marginal distribution function $F_{X|\btheta_R,\tau^2}$ can be obtained through a convolution:
\begin{equation}\label{convolution}
        1-F_{X|\btheta_R,\tau^2}(x)=P(X^*+\epsilon>x)=\int_{-\infty}^\infty \{1-F_{X^*|\btheta_R}(x-\epsilon)\}\cdot\phi_{\tau}(\epsilon) d\epsilon,
\end{equation}
where $\phi_{\tau}$ is the density of $N(0,\tau^2)$, and $F_{X^*|\btheta_R}$ is the marginal distribution function of the smooth process, whose forms are derived in \citet{huser2017bridging} and \citet{huser2017modeling} for the two models of interest in this paper. Because \eqref{convolution} cannot be further simplified, we compute the improper integral numerically using the QUADPACK algorithms which are implemented within the \texttt{gsl}\_\texttt{integration} library in \texttt{C++}. The computations in \texttt{C++} and \texttt{R} are interfaced using the package \texttt{Rcpp}. 

To compute $p$th quantile of $X$, i.e. $F_{X|\btheta_R,\tau^2}^{-1}(p)$, which is required for likelihood function evaluations, we first evaluate the distribution function $F_{X|\btheta_R,\tau^2}$ at a fine grid of $x$ values.  We then perform cubic spline interpolation through the control points to yield a continuous quantile function estimate.  Due to the smoothness of the quantile function, numerical experiments showed that this technique suffered no measurable reduction in accuracy relative to the much slower technique of computing quantiles using a numerical root finder.

\section{Additional Diagnostics}\label{sec:diagnotics}
\subsection{Marginal parameters}\label{sec:trendsurf}
\begin{figure}
    \centering
    \includegraphics[width=0.45\linewidth]{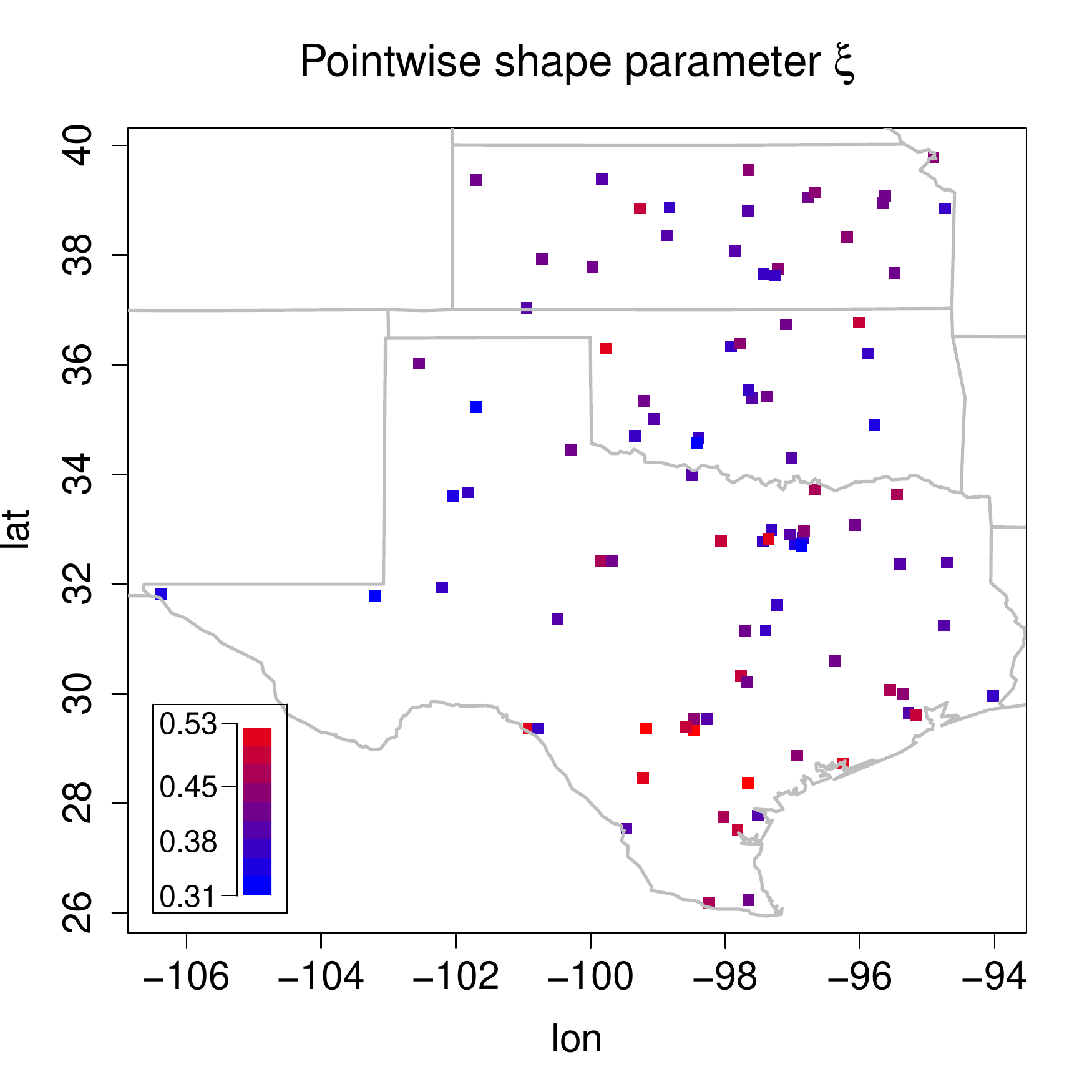}
    \includegraphics[width=0.45\linewidth]{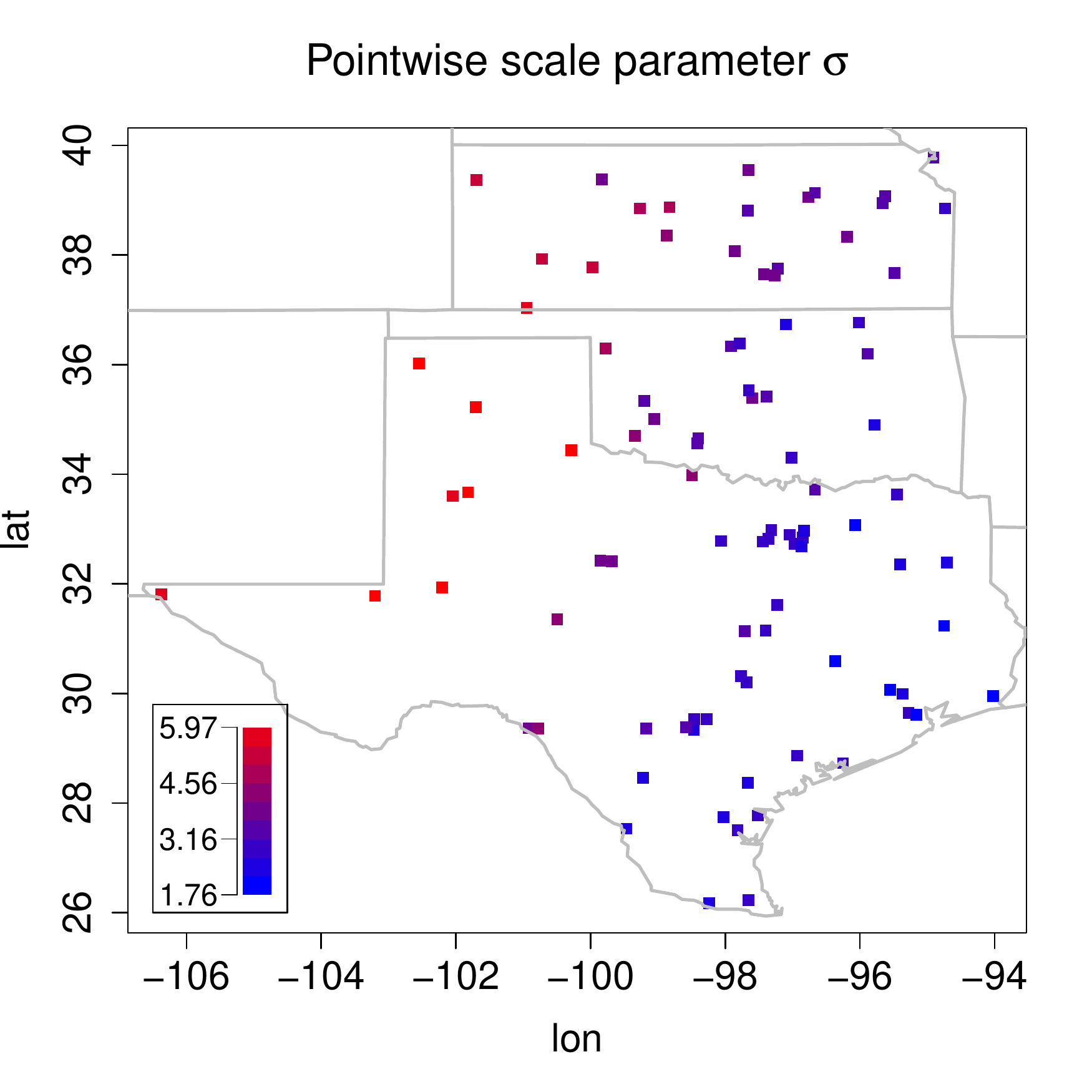}
    \caption{Pointwise generalized Pareto parameters fitted with observations from each station.}
    \label{prelim_gpd}
\end{figure}

To examine the trend surfaces of the marginal parameter $\boldsymbol{\theta}_{GPD}$, we fit univariate generalized Pareto distribution to the spring observations at each station over a high threshold $u_{0.8}$. The estimated parameters are plotted in Figure \ref{prelim_gpd}. We can see that there is no obvious spatial pattern for the shape parameter, while there is significant longitudinal effect for the scale parameter. This give grounds for applying linear trend surface to scale, and constant surface to shape (see Equation \eqref{trendsurf}).

\subsection{MCMC results}\label{mcmc_results}
\begin{figure}
    \centering
    \includegraphics[width=0.48\linewidth]{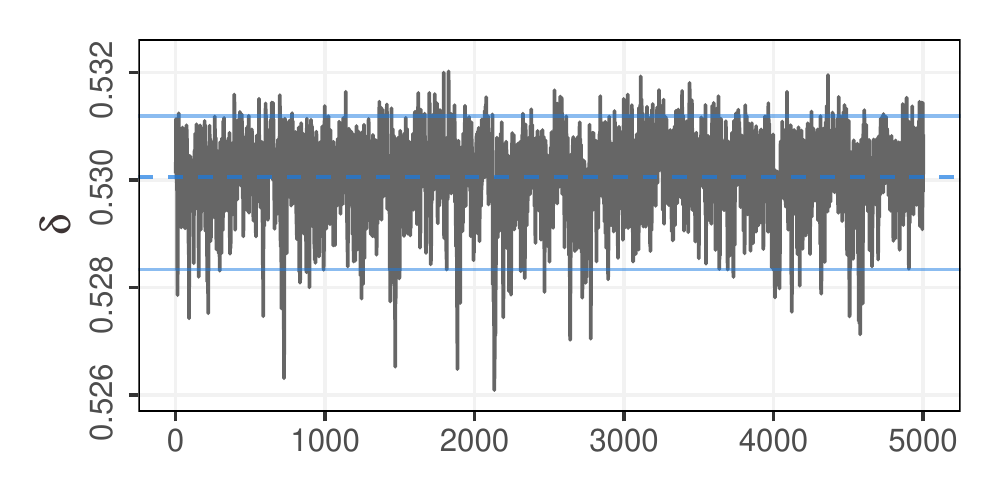}
    \includegraphics[width=0.48\linewidth]{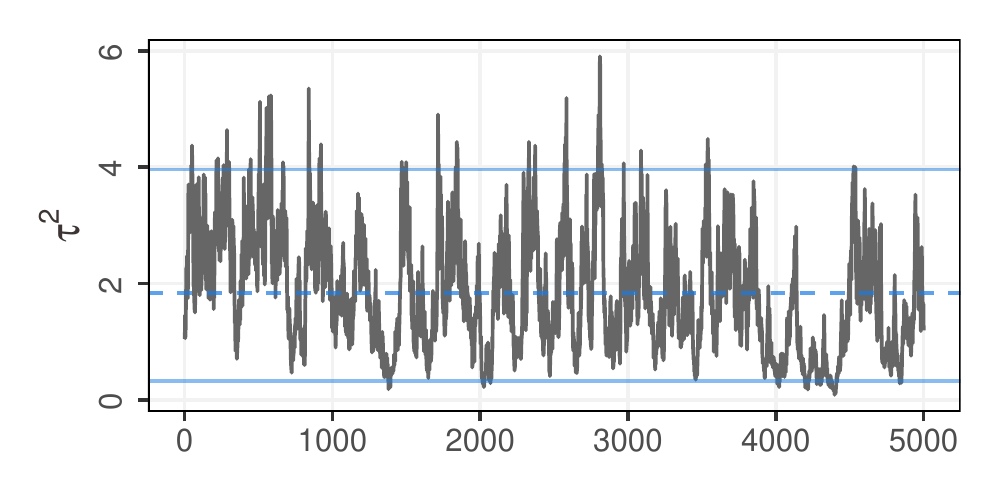}
\caption{Trace plots for $\delta$ and $\tau^2$.}
\end{figure}
\begin{figure}
    \centering
    \includegraphics[width=0.315\linewidth]{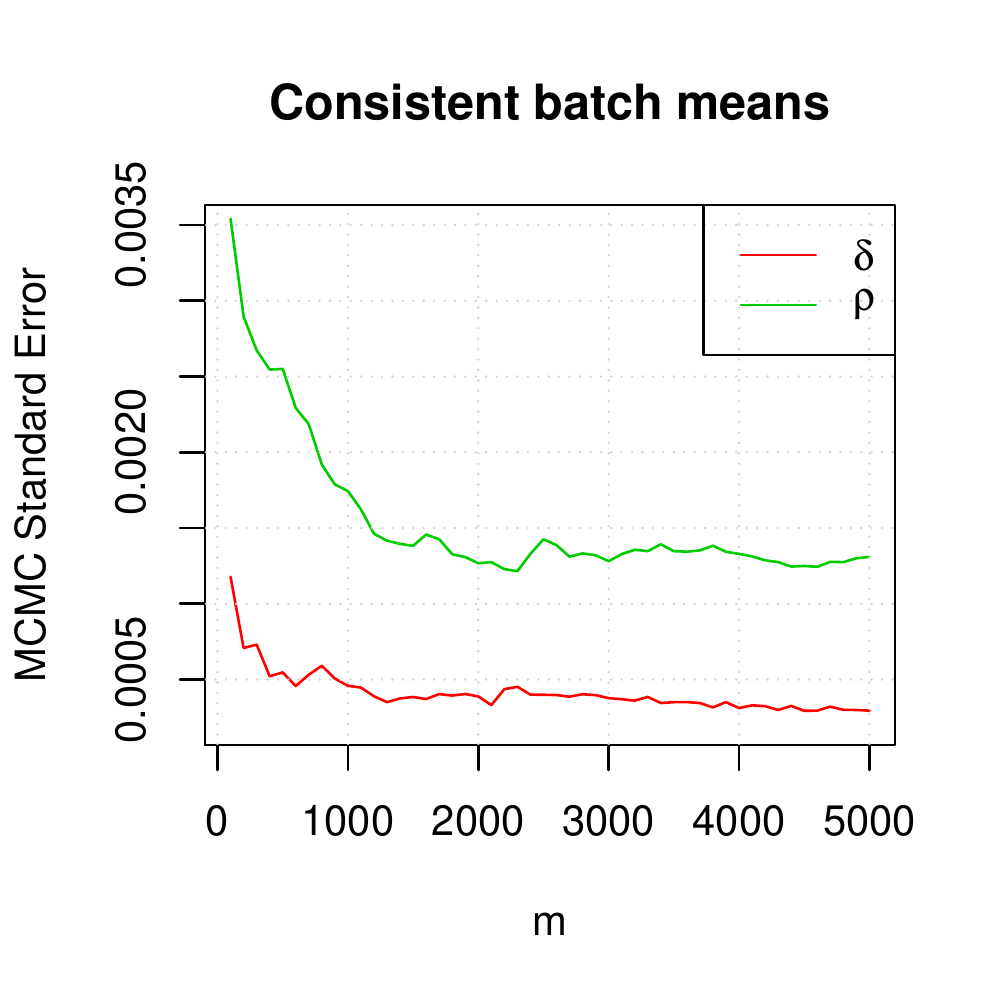}
    \includegraphics[width=0.315\linewidth]{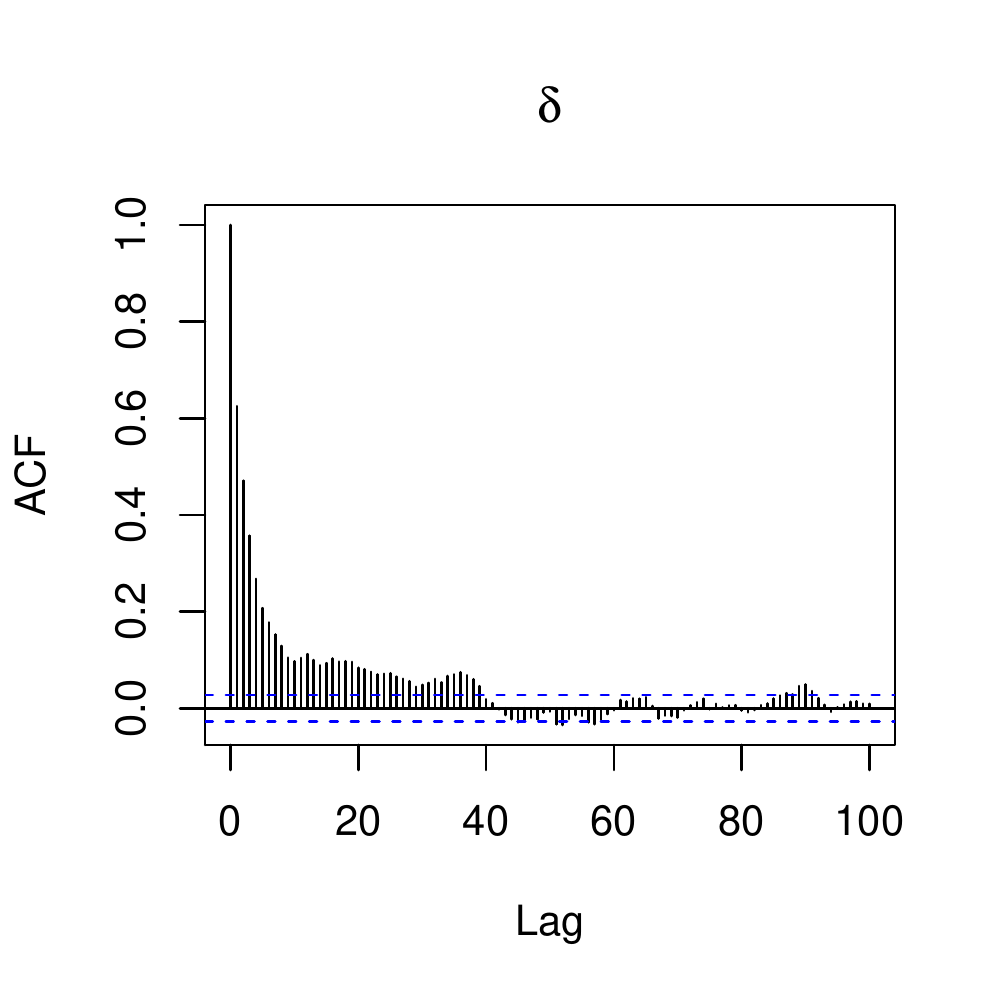}
    \includegraphics[width=0.315\linewidth]{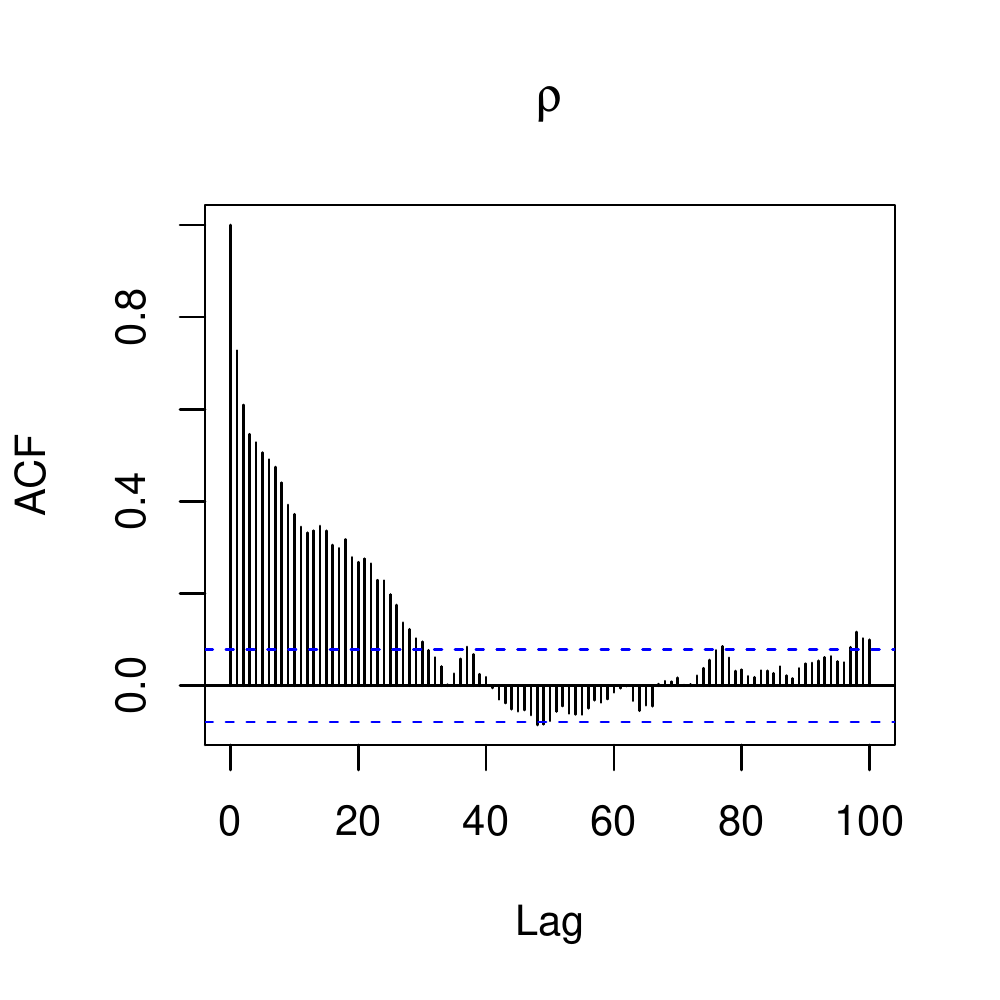}
    \caption{MCMC standard error plot (left) for $\rho$ and $\delta$, which is calculated using consistent batch means for MCMC chains. ACF plots for $\rho$ and $\delta$ are displayed on the right.}
    \label{mcmc_diag}
\end{figure}

Batch means \citep{flegal2010batch} is a convenient way to compute Monte Carlo standard errors for MCMC outputs. If one divides a Markov chain $\{X_n\}$ into $k$ batches of size $b$, the Monte Carlo estimate of $\mu=E(g(X))$ based on $i$th batch can be obtained as follows:
$$\hat{\mu}_i=\frac{\sum_{s=(i-1)b+1}^{ib}g(X_s)}{b}, i=1,\ldots, k,$$
and batch means estimate of the Monte Carlo standard error can be defined as
$$\hat{\sigma}^2=\frac{1}{k(k-1)}\sum_{i=1}^k (\hat{\mu}_i-\hat{\mu})^2,$$
where $\hat{\mu}$ is the overall Monte Carlo estimate. See Figure \ref{mcmc_diag} for batch means standard errors computed periodically for $\delta$ and $\rho$. Other parameters have similar results. We report the stabilized batch means standard errors in Table \ref{mcmc_batchMeans_sd}.

For the data analysis, we used 20 cores from the CyberScience Advanced Cyber Infrastructure at Penn State. On average, each iteration takes approximately 1.76 seconds (CPU time). The effective sample size (ESS) per second is also reported in Table \ref{mcmc_batchMeans_sd}, in which the marginal parameters has lower values.
\begin{table}
\begin{tabular}{ccccccccc}\hline
 & $\rho$ & $\nu$ & $\tau$ & $\delta$ & $\beta_0$ & $\beta_1$ & $\beta_2$ & $\xi$ \\\hline
\textbf{Monte Carlo SE} & 0.0013 & 0.0071 & 0.3570 & 0.0003 & 0.0134 & 0.0002 & 0.0001 & 0.0004\\
\textbf{ESS per second} & 0.223 & 0.468 & 0.098 & 0.190 & 0.027 & 0.027 & 0.031 & 0.019\\
\hline
\end{tabular}
\caption{Monte Carlo standard error for estimating the mean (computed using stabilized batch means for MCMC chains).}
\label{mcmc_batchMeans_sd}
\end{table}

\section*{Supplementary Material}
In this document, we examine the performance of the MCMC algorithm when $\tau^2 \rightarrow 0$, i.e., our proposed model \eqref{ourModel} converges to the underlying smooth process from \citeauthor{huser2017modeling}. One may suspect that the Markov chains will mix very poorly 
when the value of $\tau^2$ is small (although when $\tau=0$, an alternative MCMC scheme might be possible, wherein components of $\mathbf{X}^*$ could be updated from truncated distributions, albeit more complicated ones). We conduct another simulation study for datasets generated from $\delta=0.5$ (transition point) and various $\tau^2$ values ($9,\; 3,\; 0.5,\; 0.0025$). Other parameter settings remain the same as in Section \ref{paramEstimate}. \edit{The } sites and the smooth processes are simulated using the same random number generating seed for each dataset. We want to see the critical point of $\tau^2$ for the algorithm to fail.

\begin{figure}
    \centering
    \begin{subfigure}[b]{\linewidth}
    \centering
    \caption{\textbf{Case} $\tau^2=9$}
    \vskip -0.3cm
    \includegraphics[width=0.48\linewidth]{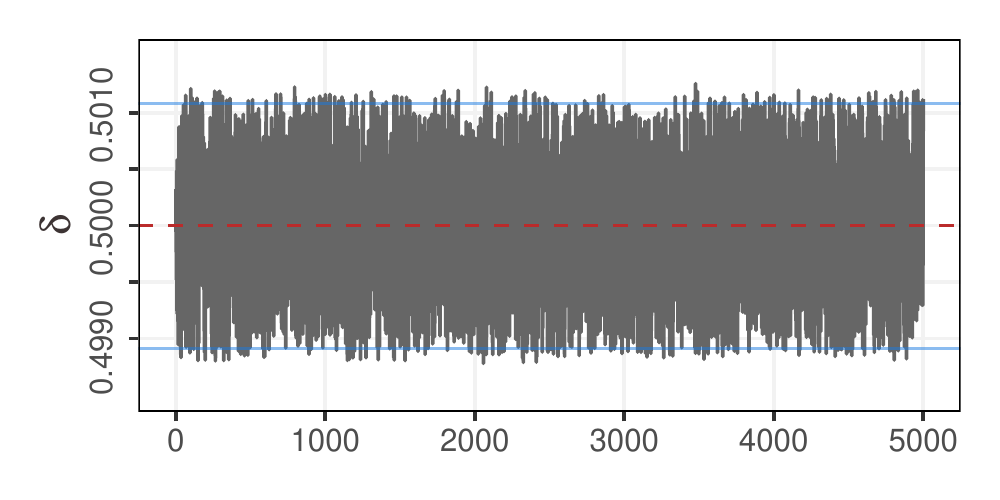}
    \includegraphics[width=0.48\linewidth]{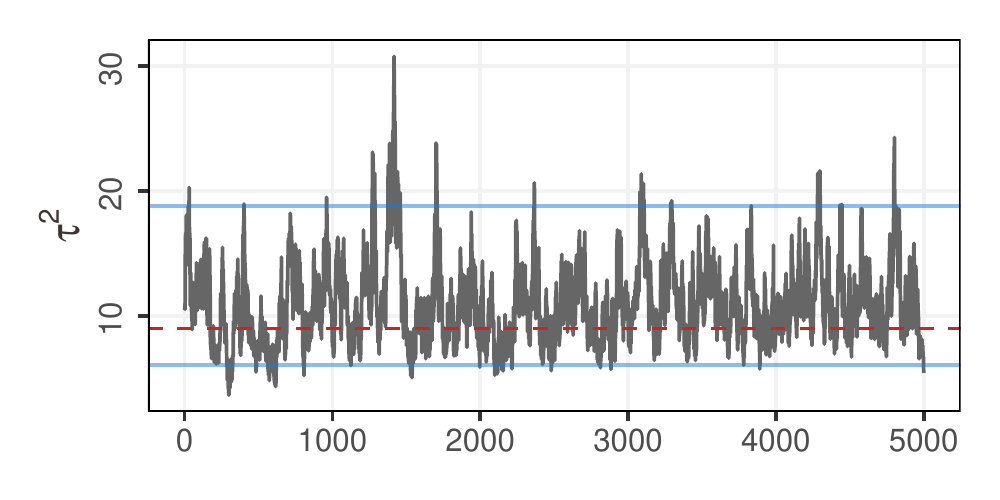}
    \end{subfigure}
    
    \begin{subfigure}[b]{\linewidth}
    \centering
    \caption{\textbf{Case} $\tau^2=3$}
    \vskip -0.3cm
    \includegraphics[width=0.48\linewidth]{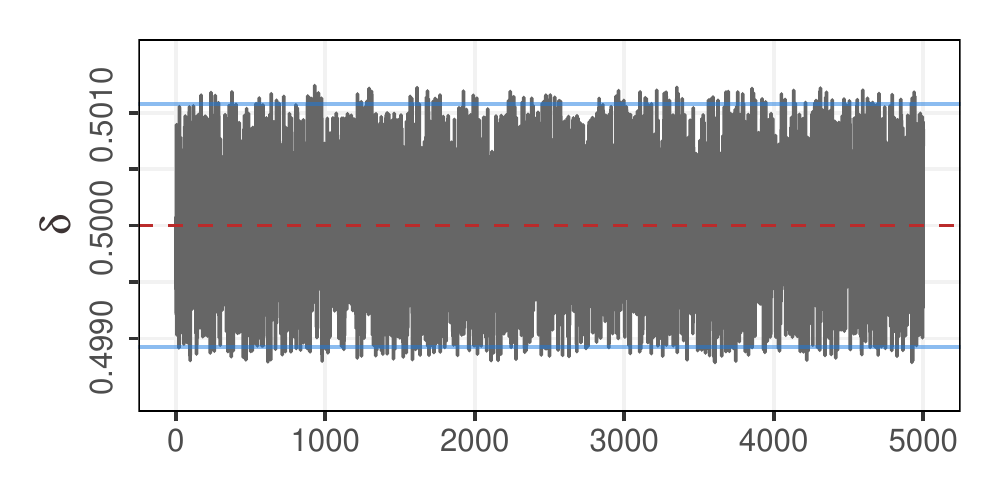}
    \includegraphics[width=0.48\linewidth]{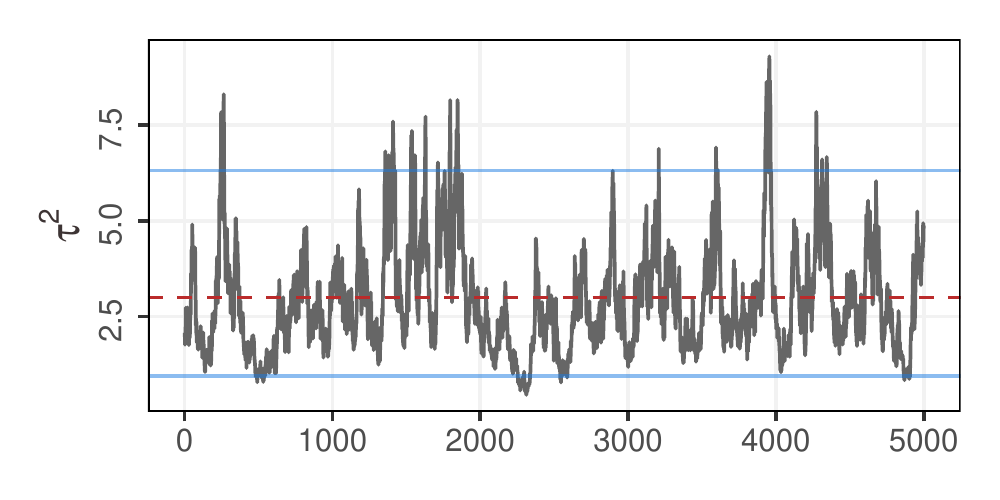}
    \end{subfigure}
    
    \begin{subfigure}[b]{\linewidth}
    \centering
    \caption{\textbf{Case} $\tau^2=0.5$}
    \vskip -0.3cm
    \includegraphics[width=0.48\linewidth]{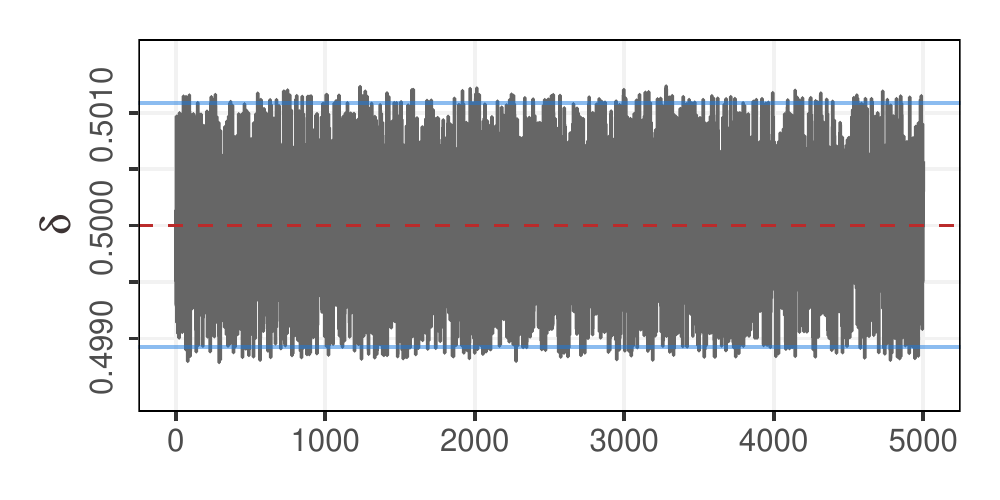}
    \includegraphics[width=0.48\linewidth]{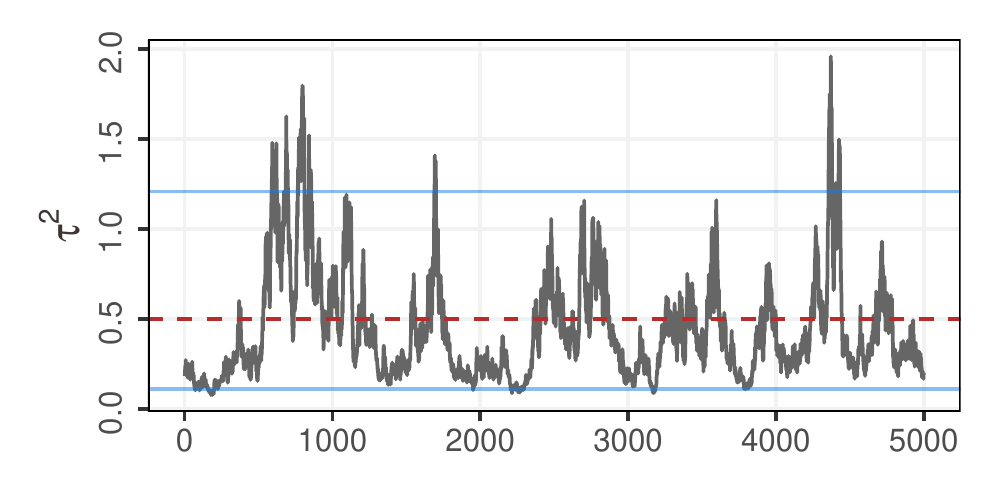}
    \end{subfigure}
    
    \begin{subfigure}[b]{\linewidth}
    \centering
    \caption{\textbf{Case} $\tau^2=0.0025$}
    \vskip -0.3cm
    \includegraphics[width=0.48\linewidth]{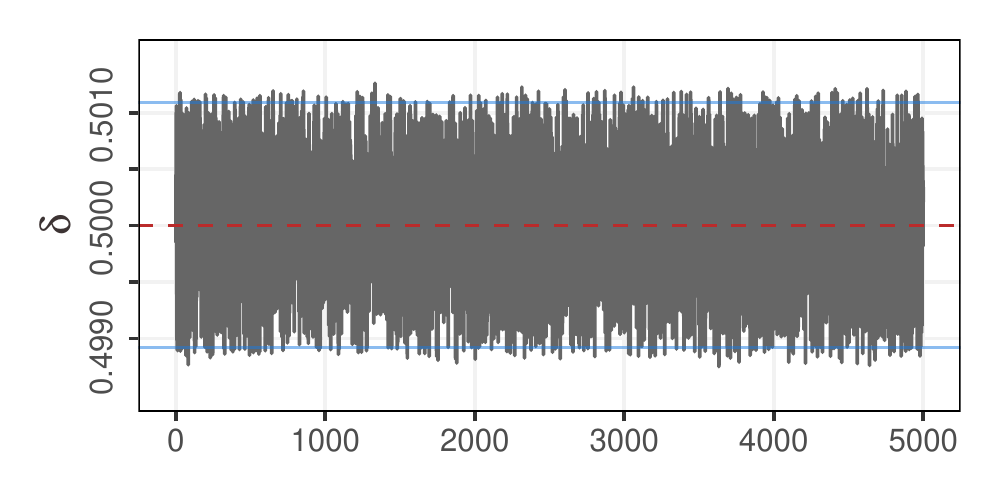}
    \includegraphics[width=0.48\linewidth]{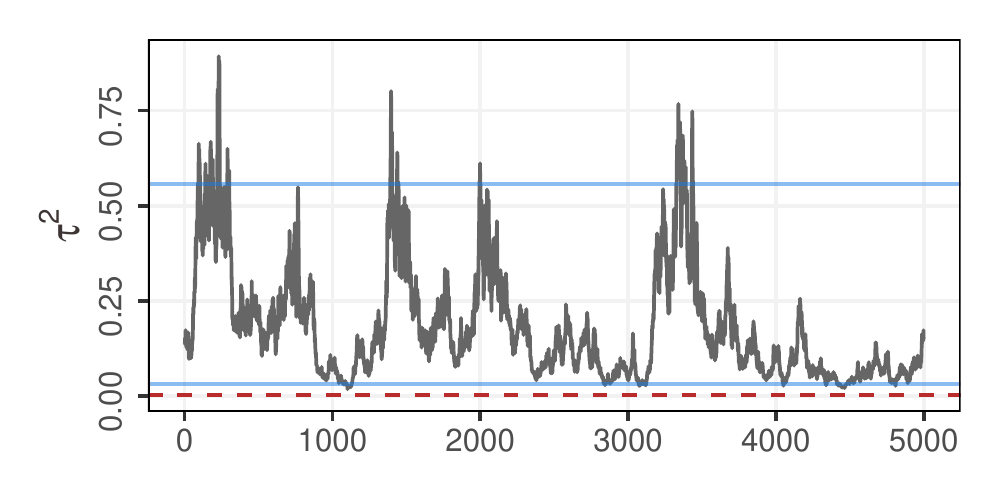}
    \end{subfigure}
    \vskip -0.4cm
    \caption{Comparisons of the trace plots of $\delta$ and $\tau^2$. See text for explanations.}
    \label{SuppMaterial_effect_of_tau}
\end{figure}
\begin{figure}
    \centering
    \begin{subfigure}[b]{0.48\linewidth}
    \centering
    \caption{\textbf{Case} $\tau^2=9$}
    \vskip -0.3cm
    \includegraphics[width=\linewidth]{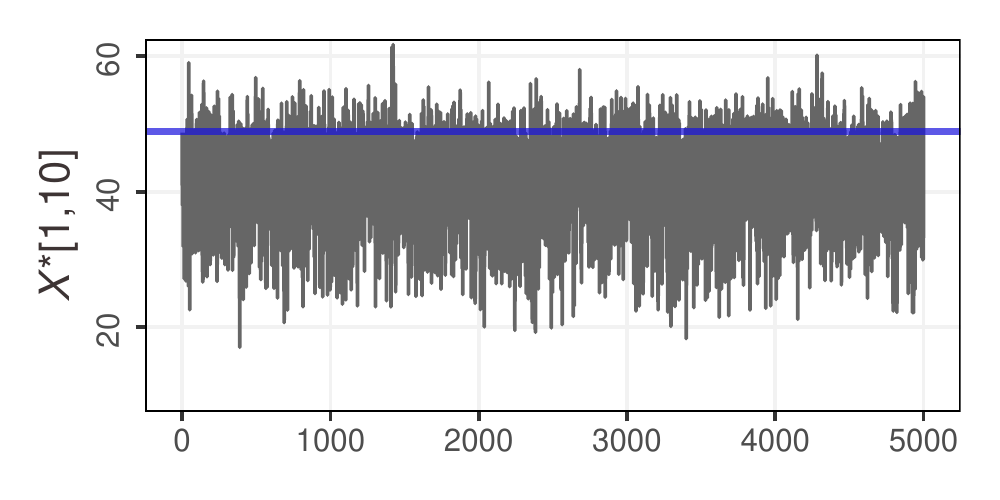}
    \end{subfigure}
    \begin{subfigure}[b]{0.48\linewidth}
    \centering
    \caption{\textbf{Case} $\tau^2=3$}
    \vskip -0.3cm
    \includegraphics[width=\linewidth]{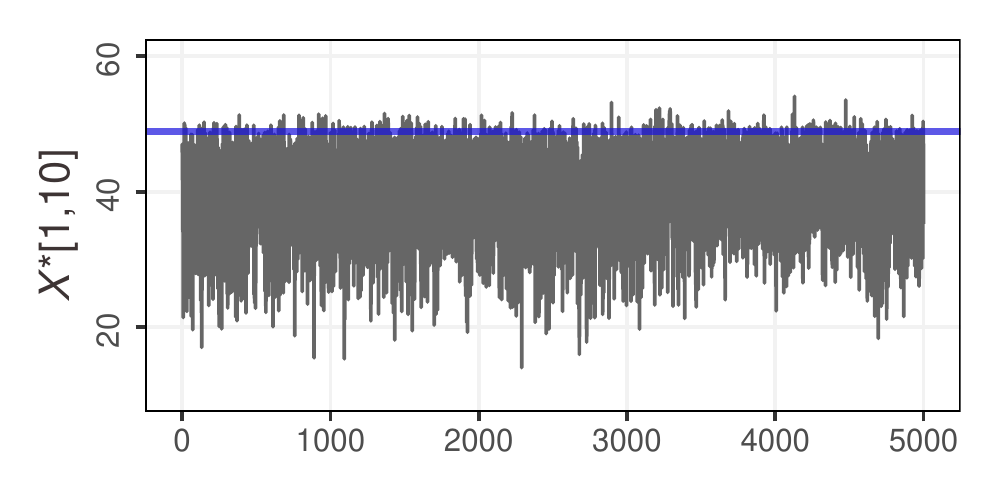}
    \end{subfigure}
    
    \begin{subfigure}[b]{0.48\linewidth}
    \centering
    \caption{\textbf{Case} $\tau^2=0.5$}
    \vskip -0.3cm
    \includegraphics[width=\linewidth]{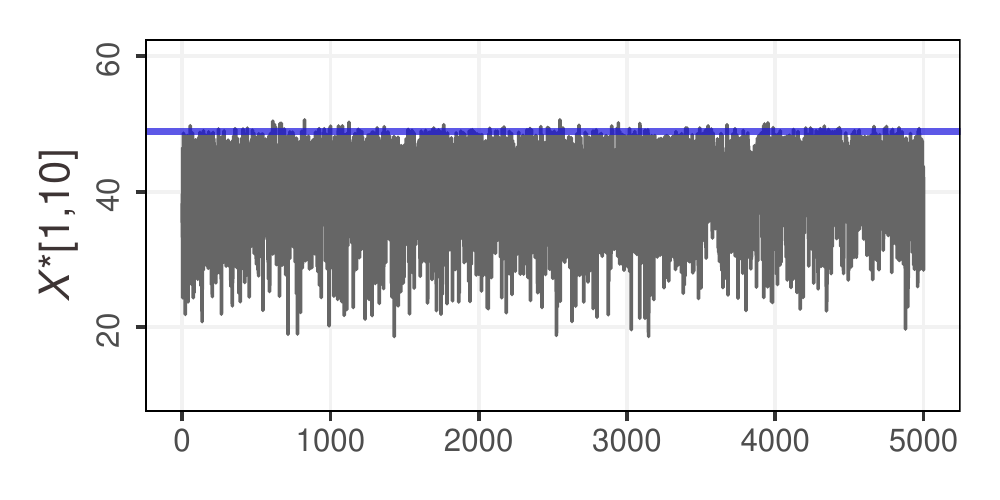}
    \end{subfigure}
    \begin{subfigure}[b]{0.48\linewidth}
    \centering
    \caption{\textbf{Case} $\tau^2=0.0025$}
    \vskip -0.3cm
    \includegraphics[width=\linewidth]{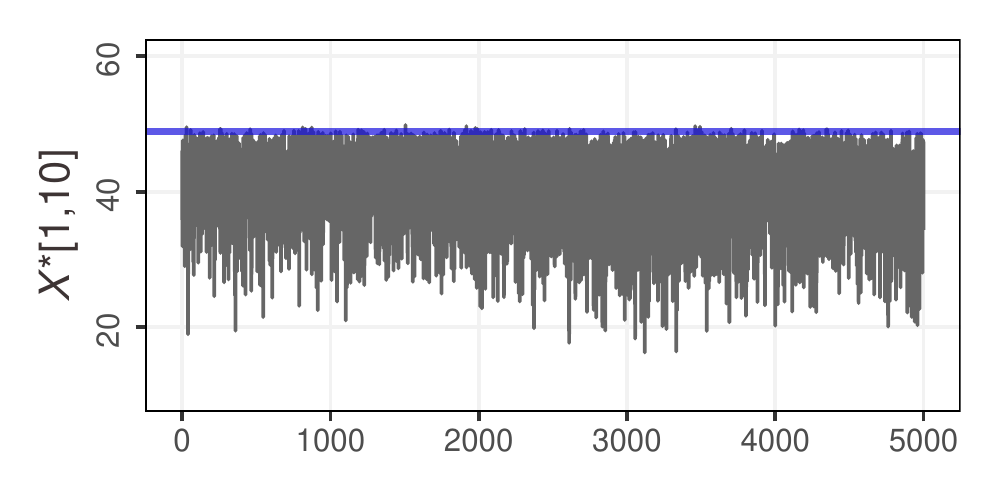}
    \end{subfigure}
    \vskip -0.4cm
    \caption{Comparisons of the trace plots of $X^*$ at 1st site and 10th time replicate. The sites and the smooth processes are the same across four simulations, and true $X^*[1,10]$ is censored.}
    \label{SuppMaterial_effect_of_tau_Xs}
\end{figure}

Figure \ref{SuppMaterial_effect_of_tau} displays the results from running the MCMC algorithm, with each row showing one case. The red dashed line signifies the true parameter values, and the blue lines indicates the 95\% posterior credible intervals. Each MCMC chain was run for 400,000 iterations. Thinning the results by a factor of 10, we show the last 50,000 iterations. We can see that, while the performance of $\delta$ stays stable for different true $\tau^2$ values, the Markov chain for $\tau^2$ converges slower when its true values becomes smaller. In the case where $\tau^2=0.0025$, the true value is outside of the 95\% credible interval and the chain mixes very poorly. Figure \ref{SuppMaterial_effect_of_tau_Xs} shows trace plots of an $X^*$ at one specific site and time which is censored for all four simulations due to the same collection of sites and underlying smoothing processes. We can see that for larger value of $\tau^2$, there are more fluctuations for $X^*$ over the threshold, which might be the reason why the MCMC algorithm works better in this case. The dissatisfying performance when $\tau^2$ is very small might also have something to do with the particular sampler we used (random walk M-H with a symmetric proposal kernel), given the parameter is so close to the boundary, rather than a problem with the overall modeling scheme.

\end{document}